\documentclass[lettersize,journal]{IEEEtran}
\usepackage{amsmath,amssymb,amsfonts}
\usepackage{algorithmic}
\usepackage{algorithm}
\usepackage{array}
\usepackage[caption=false,font=normalsize,labelfont=sf,textfont=sf]{subfig}

\usepackage{textcomp}
\usepackage{stfloats}
\usepackage{url}
\usepackage{verbatim}
\usepackage{graphicx}
\usepackage{cite}
\usepackage[colorlinks]{hyperref}
\usepackage{cleveref}
\usepackage{bbding}
\usepackage{booktabs, bm, bbm}
\usepackage[edges]{forest}
\usepackage{enumitem}
\hypersetup{linkcolor=blue,citecolor=blue}

\usepackage{CJK}

\usepackage[normalem]{ulem}

\def\ie{\emph{i.e.}}
\def\eg{\emph{e.g.}}

\usepackage{multirow, colortbl, threeparttable, makecell}

\usepackage{amsmath}
\allowdisplaybreaks[4]

\usepackage[normalem]{ulem}

\usepackage{indentfirst}

\begin{document}
\begin{CJK*}{UTF8}{gkai}
\title{Survey on Neural Routing Solvers}

\author{
Yunpeng Ba, Xi Lin, Changliang Zhou, Ruihao Zheng, Zhenkun Wang,~\IEEEmembership{Senior Member,~IEEE},\\Xinyan Liang,~\IEEEmembership{Member,~IEEE}, Zhichao Lu,~\IEEEmembership{Member,~IEEE}, Jianyong Sun,~\IEEEmembership{Senior Member,~IEEE},\\Yuhua Qian,~\IEEEmembership{Member,~IEEE}, and Qingfu Zhang,~\IEEEmembership{Fellow,~IEEE}

\thanks{Yunpeng Ba, Changliang Zhou, Ruihao Zheng, and Zhenkun Wang are with Guangdong Provincial Key Laboratory of Fully Actuated System Control Theory and Technology, School of Automation and Intelligent Manufacturing, Southern University of Science and Technology, Shenzhen, China. (E-mail: \{bayp2024, zhoucl2022, zhengrh2024\}@mail.sustech.edu.cn, wangzk3@sustech.edu.cn).}
\thanks{Xi Lin and Jianyong Sun are with the School of Mathematics and Statistics, Xi'an Jiaotong University, Xi'an, China. (E-mail: \{xi.lin, jy.sun\}@xjtu.edu.cn).}
\thanks{Xinyan Liang and Yuhua Qian are with the Institute of Big Data Science and Industry, Shanxi University, Taiyuan, China. (E-mail: liangxinyan48@163.com, jinchengqyh@126.com).}
\thanks{Zhichao Lu and Qingfu Zhang are with the Department of Computer Science, City University of Hong Kong, Hong Kong SAR, China. (E-mail: \{zhichao.lu, qingfu.zhang\}@cityu.edu.hk).}
\thanks{Corresponding author: Zhenkun Wang.}

}



\maketitle

\begin{abstract}

Neural routing solvers (NRSs) that leverage deep learning to tackle vehicle routing problems have demonstrated notable potential for practical applications. By learning implicit heuristic rules from data, NRSs replace the handcrafted counterparts in classic heuristic frameworks, thereby reducing reliance on costly manual design and trial-and-error adjustments. This survey makes two main contributions: (1) The heuristic nature of NRSs is highlighted, and existing NRSs are reviewed from the perspective of heuristics. A hierarchical taxonomy based on heuristic principles is further introduced. (2) A generalization-focused evaluation pipeline is proposed to address limitations of the conventional pipeline. Comparative benchmarking of representative NRSs across both pipelines uncovers a series of previously unreported gaps in current research.

\end{abstract}


% \IEEEraisesectionheading{\section{Introduction}}
\section{Introduction}
\label{Section 1: Intro}

The \underline{v}ehicle \underline{r}outing \underline{p}roblem (VRP) \cite{dantzig1959truck,clarke1964scheduling} is a classic \underline{c}ombinatorial \underline{o}ptimization \underline{p}roblem (COP) that seeks cost-minimizing routes for serving geographically distributed customers under specific constraints. Its scientific significance and broad practical impact have been demonstrated across various fields~\cite{vidal2013heuristics}, such as transportation~\cite{cordeau2003tabu}, logistics~\cite{hemmelmayr2012adaptive}, and manufacturing~\cite{alegre2007optimizing}. As an NP-hard problem \cite{lenstra1981complexity}, VRPs cannot be solved to optimality in polynomial time, which has driven decades of research into heuristic algorithms to obtain high-quality approximations within acceptable computation time \cite{laporte1992vehicle}. However, designing effective heuristics requires substantial domain expertise and careful manual tuning, which poses significant challenges for real-world applications.

Efforts to automate heuristic design for combinatorial optimization have long been underway. A prominent direction is algorithm selection~\cite{rice1976algorithm}, which leverages features across problem instances to choose the most suitable algorithm for a given one. This idea has been extended to portfolio-based methods~\cite{huberman1997economics, gomes2001algorithm}, where a set of complementary algorithms is maintained and selectively applied to maximize performance for different instances~\cite{huberman1997economics, gomes2001algorithm}. Another established paradigm is algorithm configuration, which aims to optimize algorithm performance for a target problem by automatically tuning parameters and combining modules~\cite{hutter2009paramils}. Despite their advances, these approaches remain confined to manually specified components and predefined parameter ranges, thus unable to discover or integrate novel algorithmic elements, which fundamentally limits potential performance gains.

\begin{table*}[htbp]
  \centering
  \caption{Comparisons of Existing Surveys for NRSs}
  \tabcolsep=0.28cm
  % \vspace{-0.2cm}
  \begin{threeparttable}
    \begin{tabular}{p{4.5em}ccp{40em}}
    \toprule[0.5mm]
     \multicolumn{1}{p{4.5em}}{\textbf{Exisiting Surveys}} & \multicolumn{1}{p{5em}}{\textbf{Unified Perspective$^\dagger$}} &  \multicolumn{1}{p{7.1em}}{\textbf{Algorithm-Level Perspective$^\ddagger$}} & \multicolumn{1}{l}{\textbf{Description}} \\
    \midrule
    Learning\newline{}Perspective  & \checkmark & $\times$  & \textbf{\textbullet} \textbf{Model Structures}: GNN, Transformer, etc.~\cite{veres2019deep,vesselinova2020learning,peng2021graph,shahbazian2024integrating,tao2025combinatorial,alanzi2025solving};\newline{}\textbf{\textbullet} \textbf{Learning Paradigms}: SL, RL, etc.~\cite{vesselinova2020learning,bengio2021machine,cappart2023combinatorial,bogyrbayeva2024machine,jin2024unified,martins2025systematic};\newline{}\textbf{\textbullet} \textbf{Generation Paradigms}: AR, NAR, etc.~\cite{peng2021graph,jin2024unified};\newline{}\textbf{\textbullet} \textbf{Reliance on Learned Modules}: E2E, Hybrid, etc.~\cite{bengio2021machine,kotary2021end,bogyrbayeva2024machine,shahbazian2024integrating,chung2025neural,alanzi2025solving}. \\
    \midrule
     Hybrid\newline{}Perspective\newline{}(Heuristic + Others) & $\times$ & \checkmark    & \textbf{\textbullet} \textbf{Residual Categories}: Grouping NRSs with \textbf{failed-to-identify heuristic categories} as\newline{}\hspace*{0.9em}``Predict''~\cite{wu2024neural} / ``Hybrid''~\cite{liu2023good,sui2025survey} / ``Non-DRL''~\cite{zhang2023review} / ``Decomposition''~\cite{bai2023analytics};
     \newline{}\textbf{\textbullet} \textbf{Self-Contradiction}: Treating \textbf{iterative} GNN-based NRSs as \textbf{construction-based}~\cite{li2022overview,wang2024solving,zhou2025learning};\newline{}\textbf{\textbullet} \textbf{Scope Limitations}: Only concerning \textbf{RL}-based~\cite{mazyavkina2021reinforcement,wang2021deep,li2021research,zong2025deep} or \textbf{Transformer}-based~\cite{araya2026makes}  NRSs.\\
    \midrule
    Heuristic\newline{}Perspective (Ours) & \checkmark & \checkmark & \textbf{\textbullet} \textbf{Hierarchical Taxonomy}: Detailing how solutions are constructed / improved;
    \newline{}\textbf{\textbullet} \textbf{Progression Identification}: Tracing NRSs from traditional heuristics;
    \newline{}\textbf{\textbullet} \textbf{Category-Specific and Category-agnostic Insights}: Transferring from heuristics to NRSs. \\
    \bottomrule[0.5mm]
    \end{tabular}%
    \begin{tablenotes}
      \footnotesize
      \item[$\dagger$] \textbf{Unified Perspective} refers to surveying the field through single-dimension attributes.
      \item[$\ddagger$] \textbf{Algorithm-Level Perspective} refers to surveying the field through differences among algorithms as a whole.

    \end{tablenotes}
    \end{threeparttable}
    % \vspace{-0.5cm}
  \label{tab:survey_existing}%
\end{table*}%

A recent and transformative development that can address this limitation in heuristic design automation is the emergence of \underline{n}eural \underline{r}outing \underline{s}olvers (NRSs). NRSs leverage deep learning (DL) models to learn implicit heuristic rules from data~\cite{kool2018attention}, replacing their handcrafted counterparts within heuristic frameworks. Their advantages over traditional heuristics primarily lie in two aspects: (1) reducing reliance on manual design by learning from data rather than manual trial-and-error tuning~\cite{bengio2021machine}, and (2) enabling GPU-accelerated parallel computation for problem solving.

The growing literature on NRSs has been partially reviewed in several surveys, yet a significant gap remains. Surveys organized from a learning perspective~\cite{veres2019deep,vesselinova2020learning,bengio2021machine,kotary2021end,peng2021graph,cappart2023combinatorial,bogyrbayeva2024machine,shahbazian2024integrating,chung2025neural,tao2025combinatorial,alanzi2025solving,martins2025systematic,jin2024unified} typically structure the field around DL techniques related to specific components, which cannot capture algorithmic structure and behavior of NRSs. Other surveys adopting a hybrid perspective~\cite{zhang2023review,bai2023analytics,liu2023good,wu2024neural,sui2025survey,li2022overview,wang2024solving,zhou2025learning,mazyavkina2021reinforcement,wang2021deep,li2021research,zong2025deep,araya2026makes} often introduce secondary attributes to define heuristic categories or restrict attention to NRSs with such attributes, leading to (1) incomplete coverage, forcing a residual ``others'' category, (2) ambiguous classification, resulting in self-contradictory taxonomies, or (3) limited scope, omitting NRSs without chosen attributes from discussion. Overall, these surveys lack a unified algorithm-level perspective on the field. This makes it difficult to clearly characterize algorithmic commonalities and differences among NRSs, and to provide systematic insights and future directions for each category and the field. A comparative summary is provided in Table~\ref{tab:survey_existing}.

This survey of NRSs makes two main contributions: (1) a unified algorithm-level review from the perspective of heuristics, and (2) a generalization-focused evaluation pipeline.

\textbf{Unified Algorithm-Level Review } NRSs are inherently heuristic algorithms powered by DL models. Building on this understanding, a hierarchical taxonomy of NRSs is proposed from the perspective of heuristics, organized by how NRSs construct or improve solutions. This perspective clarifies the relationships among NRSs and highlights their progression from traditional heuristics. Furthermore, both category-specific and category-agnostic insights from heuristics are introduced into NRSs to guide further development.

\textbf{Generalization-Focused Evaluation Pipeline } A new evaluation pipeline is proposed to address limitations of the conventional one, and representative NRSs are benchmarked under both pipelines. The key focus of the proposed pipeline is zero-shot in-problem generalization, a critical indicator of current progress in NRSs. Results under the new evaluation pipeline reveal that many NRSs are outperformed by simple construction-based heuristics such as nearest neighbor and random insertion, indicating that the conventional evaluation pipeline tends to be overly optimistic. In addition, two major challenges in the field are discussed, \ie, in-problem and cross-problem generalization, and related suggestions are provided.

The rest of this survey is organized as follows:
\begin{itemize}

\item Section \ref{Section 2: Preliminaries}  introduces essential preliminaries and details of the proposed NRS taxonomy.

\item Section \ref{Section 3: Construction}, \ref{Section 4: Single-solution}, and \ref{Section 5: Population-based} respectively review and analyze the categories of NRSs identified in the taxonomy.

\item Section~\ref{Section 6: Experiments} proposes a new evaluation pipeline and benchmarks representative NRSs in terms of zero-shot in-problem generalization, where the results reveal previously unreported limitations in NRS development.

\item Section \ref{Section 7: Challenges} outlines key research challenges and provides corresponding suggestions for future work.

\item Section \ref{Section 8: Conclusions} concludes this survey.

\end{itemize}
\section{Backgrounds and Proposed Taxonomy}
\label{Section 2: Preliminaries}

\subsection{The Formulation of Vehicle Routing Problem}

This subsection introduces a three-index vehicle flow formulation~\cite{toth2002vehicle} for the \underline{C}apacitated VRP (CVRP), which is extendable to different VRP variants. The formulation is defined on a complete directed graph $\mathcal{G}=(\mathcal{V},\mathcal{A})$, where the node set $\mathcal{V}=\{0,1,\ldots,n\}$ includes a depot (node $0$) and $n$ customers. Each arc $(i,j)\in\mathcal{A}$ is associated with a travel cost coefficient $c_{ij}>0$. The fleet consists of $K$ vehicles, each with a homogeneous capacity $C$. Each customer $i \in V \setminus \{0\}$ has a demand $d_i > 0$, while the depot has $d_0 = 0$. To capture the routing decisions, two sets of binary variables are used: $x_{ij}^k$ indicates whether vehicle $k$ traverses arc $(i,j)$, and $y_i^k$ indicates whether customer $i$ is served by vehicle $k$. The problem formulation can be defined as follows:

\begin{alignat}{2}
\label{eqn:f1} \mbox{min}\quad & \sum_{i \in V} \sum_{j \in V} c_{i j} \sum_{k=1}^{K} x_{i j}^{k}\\
\label{eqn:c1} \mbox
{s.t.}\quad & \sum_{k=1}^{K} y_{i}^{k}=1 \quad \forall i \in V \backslash\{0\},\\
\label{eqn:c2} & \sum_{k=1}^{K} y_{0}^{k}=K,\\
\label{eqn:c3} & \sum_{j \in V} x_{i j}^{k} = \!\sum_{j \in V} x_{j i}^{k} = \!y_{i}^{k} \quad \forall i \!\in\! V, k\!\in\!\{1,\ldots,K\},\\
\label{eqn:c4} & \sum_{i \in V} d_{i} y_{i}^{k} \leq C \quad \forall k\in\{1,\ldots,K\},\\
\label{eqn:c5} & \sum_{i \in S} \sum_{j \notin S} x_{i j}^{k} \geq y_{h}^{k} \quad \forall S \subseteq V \backslash\{0\}, h \in S, k\in\{1,\ldots,K\}, \\
\label{eqn:c6} & y_{i}^{k} \in\{0,1\} \quad \forall i \in V, k\in\{1,\ldots,K\}, \\
\label{eqn:c7} & x_{i j}^{k} \in\{0,1\} \quad \forall i, j \in V, k\in\{1,\ldots,K\}.
\end{alignat}

In this formulation, objective \eqref{eqn:f1} minimizes the total travel cost. Constraint \eqref{eqn:c1} ensures that each customer is visited exactly once. Constraint \eqref{eqn:c2} requires all $K$ vehicles to depart from the depot, and constraint \eqref{eqn:c3} enforces that a vehicle must arrive at and depart from the same customer. Constraints \eqref{eqn:c4} and \eqref{eqn:c5} impose capacity limit and route connectivity for each vehicle \( k \), respectively. Finally, constraints \eqref{eqn:c6} and \eqref{eqn:c7} specify the binary nature of the decision variables. This formulation explicitly identifies vehicle-arc assignments, facilitating the incorporation of additional constraints (\eg, time windows~\cite{toth2002vehicle}) and accommodating asymmetric cases. For undirected graphs, directed arc variables \(x_{ij}^k\) can be replaced with edge variables \(x_e^k\), where \(e \in E\) denotes an undirected edge.

\subsection{Categories of Heuristics}
\label{Section2.2: Heuristic Categories}

Traditional heuristics for solving VRPs can be primarily classified into two categories: construction-based methods and improvement-based methods~\cite{cook1997combinatorial}. 

\subsubsection{Construction-based Methods}

As shown in \figurename~\ref{fig:two_main_categories}, construction-based methods generate a complete solution from scratch. They can be further divided by whether solutions are generated directly on the original graph or on decomposed subgraphs. Specifically, \textbf{single-stage methods}~\cite{vidal2013heuristics} work on the original graph using simple strategies such as nearest neighbor, insertion~\cite{renaud2000heuristic}, and sweep algorithm~\cite{gillett1974heuristic}. In contrast, \textbf{two-stage methods}~\cite{toth2002vehicle} decompose the problem into different stages, typically separating customer assignment to vehicles from node sequencing within each route.

\subsubsection{Improvement-based Methods}

\begin{figure*}
    \centering
    \includegraphics[width=1\linewidth]{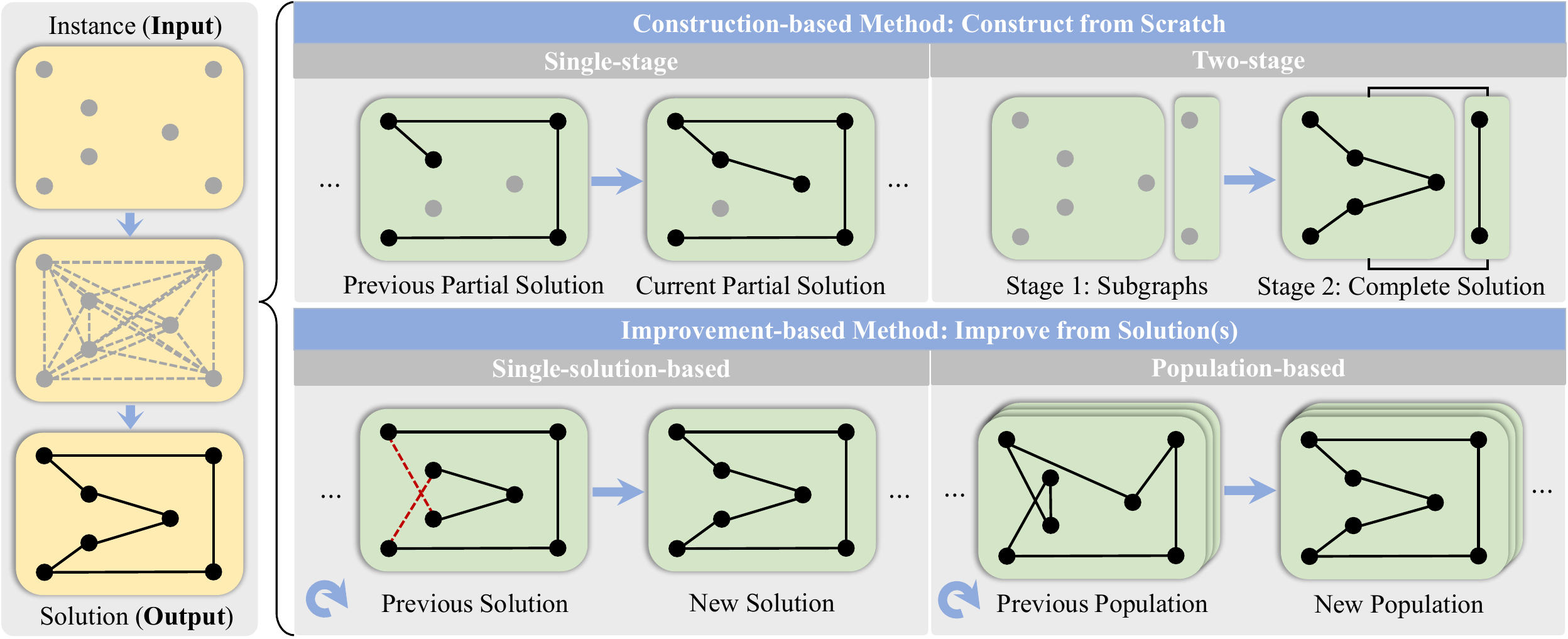}
    \caption{\textbf{Hierarchical structure of the heuristic taxonomy.} There are two main heuristic categories: the non-iterative construction-based methods and the iterative improvement-based methods. Construction-based methods can be further divided into single-stage and two-stage methods, based on whether the original graph is decomposed. Improvement-based methods can be further split into single-solution-based and population-based methods, depending on the number of solutions maintained during the improvement process.} 
    % \vspace{-0.6cm} 
    \label{fig:two_main_categories}
\end{figure*}

Improvement-based methods, also illustrated in \figurename~\ref{fig:two_main_categories},  iteratively refine one or more complete solutions during the optimization process. They can be further divided based on the number of solutions involved during the search process. Specifically, \textbf{single-solution-based methods}~\cite{funke2005local,vidal2013heuristics} focus on refining one solution by exploring its neighborhood with a small or large size. In contrast, \textbf{population-based methods}~\cite{neri2011handbook,kennedy2006swarm} maintain a population of candidate solutions, and leverage collective information to guide the search towards promising regions.

\subsection{Generation Paradigms}

NRSs mainly adopt two generation paradigms to select elements (\ie, nodes or edges): the \underline{A}uto\underline{r}egressive (AR) and \underline{N}on-\underline{a}uto\underline{r}egressive (NAR) approaches.

\subsubsection{Autoregressive Paradigm}

In this paradigm, nodes or edges are generated sequentially, with each new element conditioned on previous ones. This sequential dependency mimics a step-by-step decision-making process, which tends to yield high-quality solutions but at the cost of slower inference speed. In NRSs, the AR paradigm is well-suited to both construction-based methods that incrementally append or insert a node to a partial solution~\cite{kool2018attention}, and improvement-based methods that iteratively apply a local search move to refine a solution~\cite{wu2021learning}.

\subsubsection{Non-Autoregressive Paradigm}

In contrast, the NAR paradigm generates all elements concurrently in a single forward pass. This massively parallel strategy can significantly improve computational efficiency, though it may compromise solution quality due to the simplified independence assumption among elements. In NRSs, the NAR methods typically generate a probability distribution represented as a heatmap, over all candidate edges in the solution~\cite{joshi2019efficient}. The final solution is then generated through guided stepwise edge selection, and additional refinement steps may be applied afterwards.

\subsection{Major Learning Paradigms}

To acquire heuristic rules from data, NRSs primarily rely on two learning paradigms: \underline{S}upervised \underline{L}earning (SL) and \underline{R}einforcement \underline{L}earning (RL).

\subsubsection{Supervised Learning}

SL trains a model on a dataset of input-label pairs to learn the mapping~\cite{vinyals2015pointer,luo2023neural,drakulic2024bq}. In NRSs, the typical goal is to imitate decisions made by an expert solver, such as predicting the next node to add using a known optimal solution as the label. While this approach enables efficient learning from high-quality data, its performance is inherently bounded by the quality of the labeled data and is unlikely to surpass the expert solver it imitates.

\subsubsection{Reinforcement Learning}

RL formulates solving VRPs as a sequential decision-making process~\cite{kool2018attention,kwon2020pomo,zhou2024instance}. A solver agent learns a policy to get a high-quality solution by taking actions (\eg, selecting nodes) based on a given state (\eg, the current partial solution) and receiving a scalar reward (\eg, the negative tour length) as feedback at the end. The objective is to learn a policy that maximizes the cumulative reward. This paradigm is well-suited for NRSs as it does not require pre-solved instances, allowing the agent to explore and potentially surpass any known strategy. Nevertheless, RL-based methods may suffer from issues such as sparse rewards and high memory overhead from storing full trajectories.

\subsection{Proposed Taxonomy and Statistics}

\begin{figure*}[t]
    \centering

    \includegraphics[width=1\linewidth]{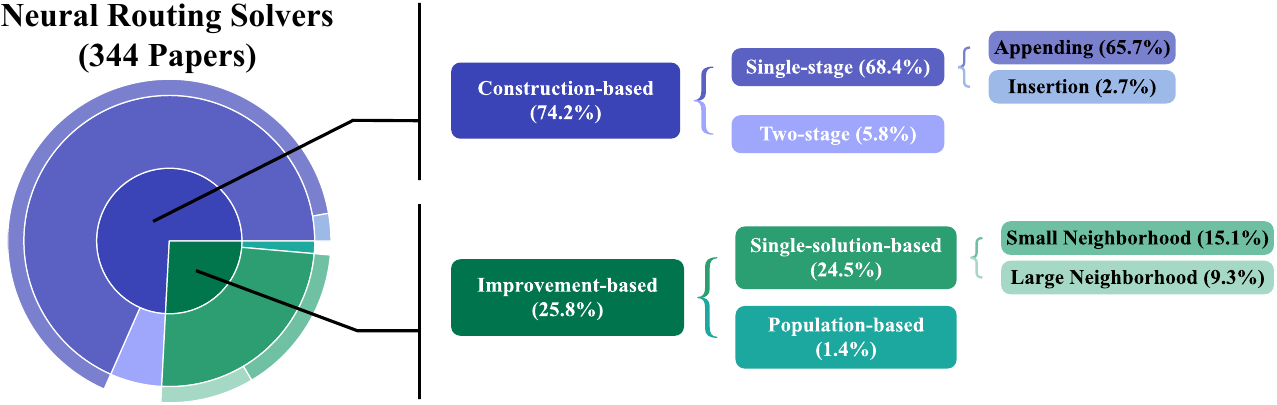}

    \caption{\textbf{Hierarchical structure of the proposed NRS taxonomy.} Each subcategory is presented with the proportion of existing studies. The statistics are obtained from Google Scholar between January 1, 2015, and November 17, 2025. The paper list is further filtered by content relevance and supplemented with relevant experience and subsequent publications. Finally, there are 439 papers, including 344 methods across categories, as well as other related studies such as surveys and benchmarks. \textbf{Note that an NRS may be counted in multiple subcategories due to its use of multiple inference strategies.}
    }
    
    % \vspace{-0.5cm}
    \label{fig:classification}
\end{figure*}

This survey proposes a taxonomy of NRSs from the perspective of heuristics. As shown in \figurename~\ref{fig:classification}, this taxonomy classifies NRSs into a multi-level hierarchy based on solution construction or improvement strategies rooted in the classical heuristic taxonomy. It naturally accommodates NRSs with different generation paradigms (AR or NAR) and learning paradigms (SL or RL) across categories. \figurename~\ref{fig:classification} further reports the statistics of all 344 NRSs across hierarchical levels.

For NRSs, the heuristic taxonomy structure outlined in Section~\ref{Section2.2: Heuristic Categories} admits finer distinctions within several subcategories. For example, single-stage NRSs for construction can be split into appending and insertion variants, depending on how nodes or edges are incorporated into a partial solution. Similarly, single-solution NRSs for improvement can be categorized by neighborhood size into small and large neighborhood methods, where the latter aligns with the traditional \underline{L}arge \underline{N}eighborhood \underline{S}earch (LNS) heuristics. Details of different categories are provided in Section \ref{Section 3: Construction}, \ref{Section 4: Single-solution}, and \ref{Section 5: Population-based}.  

\section{Construction-based Methods}
\label{Section 3: Construction}

In NRSs, construction-based methods build solutions incrementally from scratch. Similar to classical construction-based heuristics, they can be further split into single-stage and two-stage methods. The key distinction is whether the (partial) solutions are constructed directly on the original graphs or on subgraphs created by a separate decomposition stage.

\begin{figure}[t]
    \captionsetup[subfloat]{
      font=scriptsize,
      labelfont={rm,bf},
      textfont={rm,bf}
    }

    \centering
    \subfloat[Appending]{
    \begin{minipage}[b]{0.99\linewidth}
    \includegraphics[width=1\linewidth]{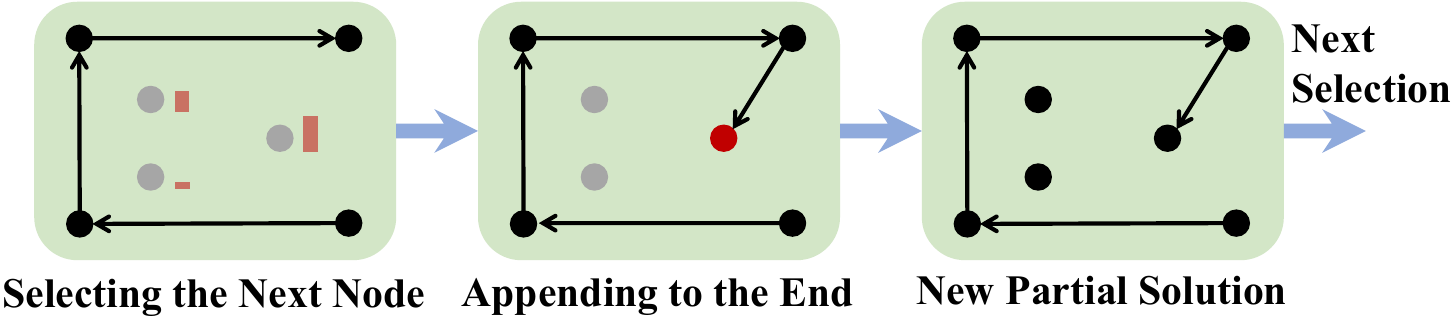}
    \end{minipage}}

    % \vspace{-0.2cm}

    \centering
    \subfloat[Insertion]{
    \begin{minipage}[b]{0.99\linewidth}
    \includegraphics[width=1\linewidth]{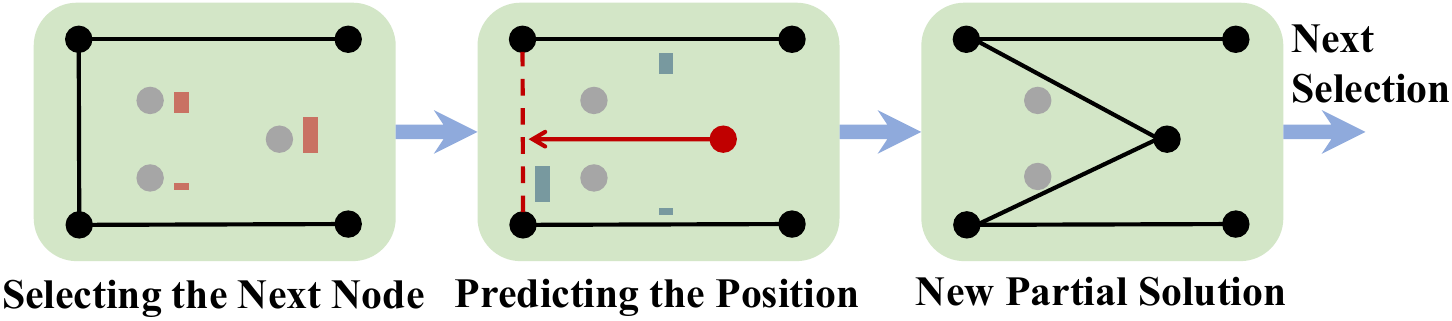}
    \end{minipage}}

    % \vspace{-0.1cm}
    \caption{\textbf{Illustration of subcategories of single-stage methods.} (a) For appending methods, the selected nodes are linked to the end of partial solutions one at a time. (b) For insertion methods, the positions are not limited.}
    % \vspace{-0.6cm}
    \label{fig:single_stage}
\end{figure}

\subsection{Single-Stage Methods}
\label{Section 3.1: Single-Stage}

Single-stage methods generate complete solutions from scratch without problem decomposition, with representative methods presented in Table~\ref{tab:single-stage}. Currently, most methods employ an appending strategy, sequentially adding selected elements to the end of the partial solution. However, alternative popular construction strategies in traditional heuristics, such as insertion, remain largely unexplored in NRSs. Therefore, although single-stage methods constitute the most active direction in current NRS research according to \figurename~\ref{fig:classification}, their design space and potential have not yet been fully explored.

\subsubsection{Appending}

As illustrated in \figurename~\ref{fig:single_stage}, nodes or edges are sequentially attached to the end of a partial solution in appending methods. Related inference strategies include greedy appending, sampling~\cite{bello2016neural,qiu2022dimes}, beam search~\cite{nazari2018reinforcement,joshi2019efficient,choo2022simulation}, (restricted) \underline{d}ynamic \underline{p}rogramming (DP)~\cite{kool2022deep}, and \underline{M}onte \underline{C}arlo \underline{t}ree \underline{s}earch (MCTS) based appending~\cite{xing2020solve}. Given that the appending position is predetermined, the core of these methods lies in stepwise element selection. Traditional heuristics such as nearest neighbor and sweep algorithms~\cite{gillett1974heuristic} rely on greedy rules based on Cartesian distance or polar angle. Corresponding NRSs replace them with learned ones.

\begin{table*}[htbp]
  \centering
  \caption{Representative Construction-based Single-stage NRSs}
  % \vspace{-0.2cm}
  \tabcolsep=0.119cm
  \begin{threeparttable}
    \begin{tabular}{llllllll}
    \toprule[0.5mm]
    \textbf{Tertiary} & \textbf{Generation} & \textbf{Solvable} & \multirow{2}[1]{*}{\textbf{Backbone}} & \textbf{Learning} & \multirow{2}[1]{*}{\textbf{Method}} & \multirow{2}[1]{*}{\textbf{Year}} & \multirow{2}[1]{*}{\textbf{Remarks}} \\
    \textbf{Category} & \textbf{Paradigm} & \textbf{VRPs} &       & \textbf{Paradigm} &       &       &  \\
    \midrule
    \multirow{37}[34]{*}{Appending} & \multirow{34}[28]{*}{AR} & \multirow{2}[2]{*}{TSP} & \multirow{2}[2]{*}{LSTM} & SL    & Ptr-Net~\cite{vinyals2015pointer} & 2015  & \multicolumn{1}{p{21.25em}}{The first NRS.} \\
\cmidrule{5-8}          &       &       &       & RL    & PN-RL~\cite{bello2016neural} & 2017  & \multicolumn{1}{p{21.25em}}{The first RL-based NRS.} \\
\cmidrule{3-8}          &       & \multicolumn{1}{p{8.5em}}{TSP, CVRP,} & \multirow{3}[0]{*}{Transformer} & \multirow{3}[0]{*}{RL} & AM~\cite{kool2018attention}    & 2019  & \multicolumn{1}{p{21.25em}}{The first NRS with the Transformer} \\
          &       & \multicolumn{1}{p{8.5em}}{OP, SDVRP,} &       &       &       &       & \multicolumn{1}{p{21.25em}}{encoder-decoder model.} \\
          &       & \multicolumn{1}{p{8.5em}}{(S)PCTSP} &       &       & MDAM~\cite{xin2021multi}  & 2021  & A multi-decoder framework with re-embedding. \\
\cmidrule{3-8}          &       & \multicolumn{1}{p{8.5em}}{(A)TSP} & Transformer & RL    & MatNet~\cite{kwon2021matrix} & 2021  & Matrix Encoding Network. \\
\cmidrule{3-8}          &       & \multicolumn{1}{p{8.5em}}{TSP, CVRP,} & \multirow{2}[1]{*}{Transformer} & \multirow{2}[1]{*}{RL} & \multirow{2}[1]{*}{Sym-NCO~\cite{kim2022sym}} & \multirow{2}[1]{*}{2022} & \multirow{2}[1]{*}{Training scheme with symmetricities.} \\
          &       & \multicolumn{1}{p{8.5em}}{PCTSP, OP} &       &       &       &       &  \\
\cmidrule{3-8}          &       & \multicolumn{1}{p{8.5em}}{PDP} & Transformer & RL    & MAPDP~\cite{zong2022mapdp} & 2022  & A multi-agent RL-based NRS for PDP. \\
\cmidrule{3-8}          &       & MOTSP, & Transformer & RL    & P-MOCO~\cite{lin2022pareto} & 2022  & A multi-objective NRS \\
          &       &   MOCVRP    &       &       &       &       & \multicolumn{1}{p{21.25em}}{with a preference-conditioned model.} \\
\cmidrule{3-8}          &       & \multicolumn{1}{p{8.5em}}{(A)TSP, CVRP, OP} & Transformer & SL    & BQ~\cite{drakulic2024bq}    & 2023  & Decoder-only structure. \\
\cmidrule{3-8}          &       & \multicolumn{1}{l}{\multirow{7}[3]{*}{TSP, CVRP}} & \multirow{7}[3]{*}{Transformer} & \multirow{2}[1]{*}{SL}    & LEHD$^\dagger$~\cite{luo2023neural}  & 2023  & Light-Encoder Heavy-Decoder structure. \\
          &       &       &       &       & SIL$^\dagger$~\cite{luo2025boosting}   & 2025  & Self-improved Training. \\
\cmidrule{5-8}          &       &       &       & \multirow{5}[0]{*}{RL} & POMO~\cite{kwon2020pomo}  & 2020  & \multicolumn{1}{p{21.25em}}{Parallel multiple rollouts.} \\
          &       &       &       &       & ELG~\cite{gao2024towards}   & 2024  & Ensemble of local and global policies. \\
          &       &       &       &       & INViT~\cite{fang2024invit} & 2024  & Distance-based search space reduction. \\
          &       &       &       &       & ICAM~\cite{zhou2024instance}  & 2025  & Distance-biased Attention. \\
          &       &       &       &       & L2R~\cite{zhou2025l2r}   & 2025  & Learning-based search space reduction. \\
\cmidrule{3-8}          &       & TSP, CVRP(TW) & Transformer & RL    & PolyNet~\cite{hottung2024polynet}   & 2025  &  Complementary solution strategies for diversity. \\
\cmidrule{3-8}          &       & min-max VRPs & Transformer & RL    & DPN~\cite{zheng2024dpn}   & 2024  & Decoupling tasks in the encoder for min-max VRPs. \\
\cmidrule{3-8}          &       & (Variants of) & \multirow{4}[0]{*}{Transformer} & \multirow{4}[0]{*}{RL} & MTPOMO~\cite{liu2024multi} & 2024  & A multi-task generalizable NRS. \\
          &       & CVRP, VRPTW,  &       &       & MVMoE~\cite{zhou2024mvmoe} & 2024  & An MoE-based NRS for multi-attribute VRPs. \\
          &       & OVRP, VRPB, &       &       & CaDA~\cite{li2024cada}  & 2025  & A constraint-prompted dual-attention mechanism. \\
          &       & VRPL &       &       & ReLD~\cite{huang2025rethinking}  & 2025  & Enhancing the Light Decoder for generalization. \\

\cmidrule{3-8}          &       & (Variants of) &  &  &  &  &  \\
          &       & (A)CVRP(TW), &       &       &       &       &  \\
          &       & OVRP, PDVRP, &    Transformer   &    RL   &   URS~\cite{zhou2025urs}    &    2025   & An NRS capable of solving more than 100 VRP  \\
          &       & VRPB(P), VRPL, &       &       &       &       & variants with a single model. \\
          &       & MDVRP, PCVRP &       &       &       &       &  \\

\cmidrule{3-8}          &       & ATSP, CVRP, & \multirow{4}[0]{*}{Transformer} & \multirow{4}[0]{*}{SL} & \multirow{4}[0]{*}{GOAL~\cite{drakulic2024goal}} & \multirow{4}[0]{*}{2025} & \multirow{4}[0]{*}{A generalist NRS with a single backbone plus} \\
          &       & CVRPTW, (S)OP, &       &       &       &       & \multirow{4}[0]{*}{problem-specific adapters.} \\
          &       & PCTSP, OVRP, &       &       &       &       &  \\
          &       & \multicolumn{1}{p{8.5em}}{SDCVRP, TRP} &       &       &       &       &  \\
\cmidrule{2-8}          &  & \multirow{2}[2]{*}{TSP} & GCN   & SL    & GCN~\cite{joshi2019efficient}   & 2019  & An NAR NRS with graph ConvNet. \\
\cmidrule{4-8}          &   NAR    &       & AGNN  & RL    & DIMES$^\dagger$~\cite{qiu2022dimes} & 2022  & Differentiable parameterizations of solution spaces. \\
\cmidrule{3-8}
&       &    TSP, CVRP   & \multicolumn{1}{p{5.125em}}{GNN} & \multicolumn{1}{p{4.25em}}{GFlowNet} & AGFN~\cite{zhang2025adversarial}  & 2025  & A GFlowNet-based construction-based NRS. \\
    \midrule
    \multirow{4}[5]{*}{Insertion} & \multirow{2}[2]{*}{AR} & TSP   & GNN   & RL    & S2V-DQN~\cite{khalil2017learning} & 2017  & A GNN-based insertion NRS. \\
\cmidrule{3-8}          &       & \multicolumn{1}{p{8.5em}}{TSP, CVRP} & Transformer & SL    & L2C-Insert$^\dagger$~\cite{luo2025learning} & 2025  & An AR SL-based insertion NRS. \\
\cmidrule{2-8}          &  \multirow{2}[2]{*}{NAR} & \multirow{2}[2]{*}{TSP} & U-Net & SL    & DMPP~\cite{graikos2022diffusion}  & 2022  & An NAR NRS with image-based diffusion models. \\
\cmidrule{4-8}          &       &       & AGNN  & SL    & DIFUSCO$^\dagger$~\cite{sun2023difusco} & 2023  & An NAR NRS with graph-based diffusion models. \\
    \bottomrule[0.5mm]
    \end{tabular}%
    \begin{tablenotes}
      \footnotesize
      \item[$\dagger$] The NRS supports multiple inference strategies and currently employs a construction-based greedy inference.
    \end{tablenotes}
  \end{threeparttable}
    % \vspace{-0.1cm}
  \label{tab:single-stage}%
\end{table*}%

In AR appending methods, learned rules sequentially select the next node to append based on the current solution state. Consequently, this approach has driven efforts to improve the model's ability for state representation and reasoning based on the current state. The first NRS Ptr-Net~\cite{vinyals2015pointer} employs an attention-based pointer mechanism for stepwise node selection. PN-RL~\cite{bello2016neural} introduces RL into NRSs and adopts active search to fine-tune on individual test instances. AM~\cite{kool2018attention} incorporates the Transformer-based encoder-decoder architecture, which improves state representation ability. Subsequently, POMO~\cite{kwon2020pomo} extends this work by leveraging multiple trajectories with different starting nodes to enhance exploration.

\begin{table*}[ht]
  \centering
  \caption{Representative Construction-based Two-stage NRSs}
    \tabcolsep=0.158cm
    % \vspace{-0.2cm}
    \begin{tabular}{p{8em}lrrrrrr}
    \toprule[0.5mm]
    \multicolumn{1}{l}{\textbf{Role of the}} & \textbf{Generation} & \multicolumn{1}{l}{\textbf{Solvable}} & \multicolumn{1}{l}{\multirow{2}[0]{*}{\textbf{Backbone}}} & \multicolumn{1}{l}{\textbf{Learning}} & \multicolumn{1}{l}{\multirow{2}[0]{*}{\textbf{Method}}} & \multicolumn{1}{l}{\multirow{2}[0]{*}{\textbf{Year}}} & \multicolumn{1}{l}{\multirow{2}[0]{*}{\textbf{Remarks}}} \\
    \multicolumn{1}{l}{\textbf{First Stage}} & \textbf{Paradigm} & \multicolumn{1}{l}{\textbf{VRPs}} &       & \multicolumn{1}{l}{\textbf{Paradigm}} &       &       &  \\
    \midrule
    Scale Reduction & AR    & \multicolumn{1}{l}{TSP} & \multicolumn{1}{l}{CNN,} & \multicolumn{1}{l}{RL} & \multicolumn{1}{l}{H-TSP~\cite{pan2023htsp}} & \multicolumn{1}{l}{2023} & \multicolumn{1}{l}{A two-stage NRS capable for TSP instances with} \\
    \multicolumn{1}{r}{} &       &       &  \multicolumn{1}{l}{Transformer}    &       &       &       & \multicolumn{1}{l}{10K nodes.} \\
    \midrule
    Scale Reduction; & AR     & \multicolumn{1}{l}{CVRP} & \multicolumn{1}{l}{Transformer} & \multicolumn{1}{l}{RL} & \multicolumn{1}{l}{TAM-AM~\cite{hou2023generalize}} & 2023  & \multicolumn{1}{l}{A two-stage NRS capable for VRP instances with} \\
    \multicolumn{1}{p{9.2em}}{Constraint Handling} &       &       &       &       &       &       & \multicolumn{1}{l}{over 5K nodes.} \\
    \bottomrule[0.5mm]
    \end{tabular}%
    % \vspace{-0.6cm}
  \label{tab:two-stage}%
\end{table*}%

To mitigate interference from irrelevant information of visited nodes during stepwise selections, some methods periodically re-embed feasible nodes~\cite{peng2020deep,xin2020step}. This alternating process of re-encoding for updated embeddings and decoding for node selection inevitably incurs substantial computational cost, yet enables more accurate state representation. A more direct alternative shifts the computational burden to a stronger decoder that performs stepwise dynamic node re-embedding. As a result, the original \underline{H}eavy \underline{E}ncoder and \underline{L}ight \underline{D}ecoder (HELD) structure is replaced by a \underline{L}ight \underline{E}ncoder \underline{H}eavy \underline{D}ecoder (LEHD) or even a Decoder-only structure. Methods adopting this design, such as BQ~\cite{drakulic2024bq} and LEHD~\cite{luo2023neural}, demonstrate improved generalization ability. These two methods typically rely on SL to achieve high sample efficiency.

For VRPs with common distance-related objectives, the optimal next node to append is often close to the current partial solution. This observation has inspired two distinct strategies. The first is to restrict candidate selection at each step to a local neighborhood based on distance, which dramatically reduces the computational complexity and difficulty of node selection with little loss in solution quality. The second is to explicitly incorporate node-wise distance information into specific modules, thereby enhancing the model's ability to assess the current state. Examples of the former local policy approach include~\cite{gao2024towards,wang2024distance,fang2024invit,drakulic2024bq,goh2024hierarchical,zhou2025l2r}. Particularly, ELG~\cite{gao2024towards} introduces an auxiliary local policy on polar-coordinate features in addition to the regular global policy, while INViT~\cite{fang2024invit} aggregates multi-scale neighborhood information through nested local views. The latter distance-enhanced modeling strategy is demonstrated in~\cite{jin2023pointerformer,son2023meta,li2023learning,gao2024towards,wang2024distance,zhou2024instance,huang2025rethinking,zhou2025l2r}. Specifically, ICAM~\cite{zhou2024instance} introduces a distance-based adaptation function within the attention mechanism to better capture spatial relationships.

AR appending remains an active research area in NRSs, with a notable advantage that lies in learning a relatively simple stepwise node-selection policy. However, like traditional construction-based heuristics, these methods generate solutions from scratch, where suboptimal selections in early steps inevitably impact the quality of subsequent decisions. The solution quality can be further improved by iteratively refining solution segments using the same inference mechanism, though at the cost of increased computational overhead. Such refinement strategies fall under the restricted direct LNS subcategory, which is discussed further in Section~\ref{Section 4.2.1.2: Direct LNS (Restricted)}.

In contrast, NAR appending methods like GCN~\cite{joshi2019efficient} and DIMES~\cite{qiu2022dimes} select elements in a single pass guided by a predicted heatmap. While this approach enables faster inference, the static nature of heatmap cannot account for the influence of dynamic masked elements or the evolving partial solution, thereby gradually distorting the guidance information and leading to suboptimal performance. To mitigate this limitation, one potential direction is to develop an inference process that dynamically updates the heatmap during element selection. Another plausible direction is to adopt iterative refinement, thereby converting these methods into improvement-based approaches. In such cases, techniques such as population-based strategies or local search can be applied to refine the solutions (see more details in Section~\ref{Section 4: Single-solution} and \ref{Section 5: Population-based}).

\subsubsection{Insertion}

Within single-stage methods, insertion remains a notable yet underexplored alternative to the prevalent appending paradigm. As illustrated in \figurename~\ref{fig:single_stage}, insertion methods can place unvisited nodes into arbitrary positions of the partial solution, rather than only at the end. This flexibility introduces two coupled decisions, namely, which nodes to insert and where to insert them. While having higher time complexity, insertion can mitigate error accumulation inherent in appending by allowing corrections in subsequent steps.

Only a few studies have attempted to learn insertion policies. Among AR insertion approaches, S2V-DQN~\cite{khalil2017learning} selects nodes with the highest predicted values and inserts them at the minimum-cost positions for TSP. Besides, L2C-Insert~\cite{luo2025learning} selects unvisited nodes via a nearest-neighbor rule and learn to determine the insertion positions. A few NAR insertion methods, such as DIFUSCO~\cite{sun2023difusco}, incorporate greedy edge insertion guided by predefined priority scores as one inference strategy. These initial efforts, however, only scratch the surface of insertion. Future research could investigate the design of more effective joint policies that explicitly model node-position interactions, while balancing computational overhead with the opportunity to repair earlier suboptimal decisions.

\subsection{Two-Stage Methods}
\label{Section 3.2: Two-Stage}

Two-stage methods are designed to address different challenges separately in each stage. The widely used ``cluster-first route-second'' strategy~\cite{gillett1974heuristic,fisher1981generalized} first groups customers into feasible clusters based on constraints and then sequences nodes within each cluster (\ie, subgraph) by solving a set of smaller TSPs. It can significantly reduce the problem scale and allows the second stage to focus on sequencing. This strategy has been adopted by a few NRSs as shown in Table~\ref{tab:two-stage}. For example, TAM-AM~\cite{hou2023generalize} partitions a large-scale VRP into clusters of small-scale TSPs in the first stage, and then applies a single-stage solver such as AM for each TSP in the second stage. Other methods, such as H-TSP~\cite{pan2023htsp}, solely target the scaling challenge of TSPs by decomposition, which generate open-loop tours per cluster and then connect them to form a complete solution. By leveraging existing TSP solvers, these methods essentially transfer the core challenge of problem solving to the graph-partitioning step. Furthermore, related improvement-based methods with iterative redivision and refinement are discussed in Section~\ref{Section 4.2.1.2: Direct LNS (Restricted)}.

\begin{figure}[t]
    \captionsetup[subfloat]{
      font=scriptsize,
      labelfont={rm,bf},
      textfont={rm,bf}
    }

    \centering
    \subfloat[Small Neighborhood]{
    \begin{minipage}[b]{0.99\linewidth}
    \includegraphics[width=1\linewidth]{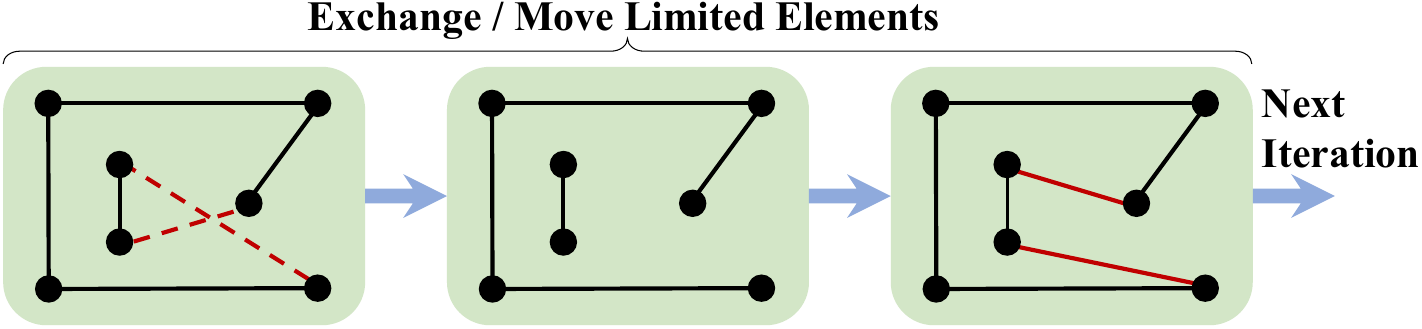}
    \end{minipage}}

    % \vspace{-0.2cm}

    \centering
    \subfloat[Large Neighborhood]{
    \begin{minipage}[b]{0.99\linewidth}
    \includegraphics[width=1\linewidth]{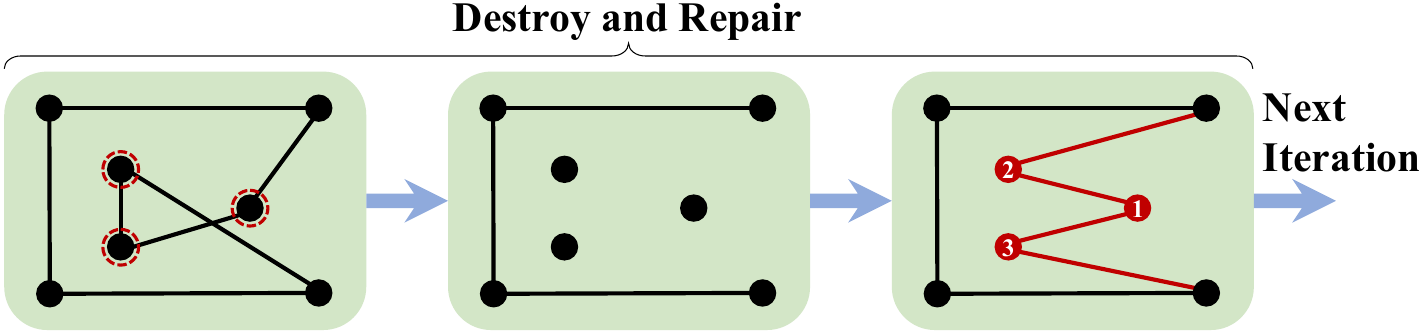}
    \end{minipage}}

    % \vspace{-0.3cm}
    
    \caption{ \textbf{Illustration of subcategories of single-solution-based methods.} (a) For small neighborhood methods, an example with the 2-opt operator is presented. At each iteration, two edges are replaced. (b) For large neighborhood methods, the classic destroy-and-repair process is illustrated. At each iteration, some nodes are picked out and then reinserted one by one in a prescribed order (indicated by the numbered labels).}
    % \vspace{-0.6cm}
    \label{fig:single_solution}
\end{figure}

\section{Single-Solution-Based Methods\\For Improvement}

\label{Section 4: Single-solution}

Single-solution-based methods iteratively improve a complete solution by exploring its neighborhood, which is a specific subset of feasible solutions reachable from the current solution through specific modifications.

\subsection{Small Neighborhood Methods}

\begin{table*}[htbp]
  \centering
  \caption{Representative Improvement-based Single-solution-based Small Neighborhood NRSs}
  \tabcolsep=0.11cm
  % \vspace{-0.2cm}
  \begin{threeparttable}
    \begin{tabular}{rrllllll}
    \toprule[0.5mm]
    \multicolumn{1}{l}{\textbf{Quaternary}} & \textbf{Generation} & \textbf{Solvable} & \multirow{2}[1]{*}{\textbf{Backbone}} & \textbf{Learning} & \multirow{2}[1]{*}{\textbf{Method}} & \multirow{2}[1]{*}{\textbf{Year}} & \multirow{2}[1]{*}{\textbf{Remarks}} \\
    \multicolumn{1}{l}{\textbf{Category}} & \multicolumn{1}{l}{\textbf{Paradigm}} & \textbf{VRPs} &       & \textbf{Paradigm} &       &       &  \\
    \midrule
    \multicolumn{1}{l}{\multirow{7}[7 ]{*}{Immediate}} & \multicolumn{1}{l}{\multirow{6}[6]{*}{AR}} & CVRP  & LSTM & RL    & NeuRewriter~\cite{chen2019learning} & 2019  & Improvement with separate policies to select node-pairs. \\
\cmidrule{3-8}          &       & \multicolumn{1}{l}{\multirow{2}[1]{*}{TSP, CVRP}} & \multicolumn{1}{l}{\multirow{2}[1]{*}{Transformer}} & \multirow{2}[1]{*}{RL} & LIH~\cite{wu2021learning}   & 2021  & Improvement with a single policy to select node-pairs. \\
          &       &       &       &       & DACT~\cite{ma2021learning}  & 2021  & Improvement with cyclic positional encoding. \\
\cmidrule{3-8}          &       & TSP   & \multicolumn{1}{p{5.565em}}{GCN, FiLM,} & RL    & Neural-3-OPT~\cite{sui2021learning} & 2021  &
 Improvement with the 3-opt operator. \\
          &       &       & \multicolumn{1}{p{5.565em}}{LSTM} &       &       &       &  \\
\cmidrule{3-8}          &       & PDP   & Transformer & RL    & NCS~\cite{kong2024efficient}   & 2024  & An improvement-based NRS for PDP. \\
\cmidrule{2-8}          & \multicolumn{1}{p{4.315em}}{NAR} & TSP   & GNN   & SL    & RGLS~\cite{hudson2022graph}  & 2022  & Predicting regret for guided local search. \\
    \midrule
          & \multicolumn{1}{l}{\multirow{2}[1]{*}{AR}} & TSP,  & Transformer,  & RL    & NeuOpt~\cite{ma2024learning} & 2023  & Improvement with flexible k-opt. \\
          &       & \multicolumn{1}{p{5em}}{CVRP} & \multicolumn{1}{p{5.565em}}{GRU} &       &       &       &  \\
\cmidrule{2-8}          &       & TSP, PDP, &       &       &       &       &  \\
          &       & CVRP, & SGN   & SL, UL & NeuroLKH~\cite{xin2021neurolkh} & 2021  & Introducing DL to LKH. \\
          &       & CVRPTW &       &       &       &       &  \\
\cmidrule{3-8}    \multicolumn{1}{l}{Sequential} &       & \multirow{6}[5]{*}{TSP} & GCRN  & SL    & \multicolumn{1}{p{6.5em}}{Att-GCN~\cite{fu2021generalize}} & 2021  & Introducing MCTS-k-opt to NAR NRSs. \\
\cmidrule{4-8}          & \multicolumn{1}{l}{NAR} &       & AGNN  & RL    & \multicolumn{1}{p{5.94em}}{DIMES$^\ddagger$~\cite{qiu2022dimes}} & 2022  & Differentiable parameterizations of solution spaces. \\
\cmidrule{4-8}          &       &       & AGNN  & SL    & \multicolumn{1}{p{5.94em}}{DIFUSCO$^\ddagger$~\cite{sun2023difusco}} & 2023  & An NAR NRS with graph-based diffusion models. \\
\cmidrule{4-8}          &       &       & SAG   & UL    & UTSP~\cite{min2024unsupervised}  & 2023  & A NAR UL-based NRS. \\
\cmidrule{4-8}          &       &       & /     & /     & SoftDist~\cite{xia2024position} & 2024  & A critique of DL-output-heatmap-MCTS-k-opt paradigm. \\
    \bottomrule[0.5mm]
    \end{tabular}%
    \begin{tablenotes}
      \footnotesize
      \item[$\ddagger$] The NRS supports multiple inference strategies and currently employs an improvement-based MCTS-k-opt inference.
    \end{tablenotes}
  \end{threeparttable}
    % \vspace{-0.6cm}
  \label{tab:small neighbor}%
\end{table*}%

As presented in \figurename~\ref{fig:single_solution}, small neighborhood methods explore neighborhoods with limited sizes defined by local search operators. Based on whether moves are decomposed, they are categorized into immediate and sequential search. Immediate search relies on simple operators such as swap and 2-opt. In contrast, sequential search employs more complex operators like k-opt (k$>$2), where each move is typically decomposed into a sequence of steps to mitigate decision complexity.

\subsubsection{Immediate Search}
\label{Section 4.1.1: Immediate}

There are typically two steps in each iteration of immediate search methods: (1) selecting a few nodes or edges, and (2) performing a single move via a local search operator. Learned rules in these methods primarily focus on the selection step, implemented either autoregressively, such as choosing nodes or edges with an agent, or non-autoregressively, such as generating a heatmap to guide iterative node pair or edge selection.

AR immediate search methods select moves via learned rules rather than handcrafted distance-based ones. For example, NeuRewriter~\cite{chen2019learning} uses two interrelated learned rules to separately select two nodes for a local search move, while LIH~\cite{wu2021learning} and DACT~\cite{ma2021learning} employ a single learned rule to select node pairs. DACT additionally addresses challenges related to positional encoding. In contrast, Neural-3-OPT~\cite{sui2021learning} learns separate rules to remove and reconnect three edges for each 3-opt move. NAR immediate search methods, such as RGLS~\cite{hudson2022graph}, use heatmaps predicted by learned rules to guide the improvement process. For example, regret values can be predicted for all edges to steer the improvement process in \underline{g}uided \underline{l}ocal \underline{s}earch (GLS)~\cite{gendreau2010handbook}.

Benefiting from fine-grained local search operators, immediate search methods typically perform well on small-scale instances. Nevertheless, their limited neighborhood size makes them prone to local optima and less effective on large-scale problems. Therefore, the development of such methods has encountered a bottleneck in recent years.

\subsubsection{Sequential Search}
\label{Section 4.1.2: Sequential}

Sequential search methods typically employ k-opt operators with k$>$2 to expand the search neighborhood for discovering better solutions. However, increasing k would lead to exponential growth in neighborhood size and, consequently, in computational complexity. A plausible strategy to address this is to decompose a k-opt move into a sequence of basic moves, which treats the improvement process as a Markov Decision Process.

There are various strategies to select a basic move at each step. To begin with, AR methods learn rules for stepwise basic move selections. For example, NeuOpt~\cite{ma2024learning} dynamically adjusts k to balance coarse- and fine-grained search. Besides, NAR methods prioritize basic moves based on per-edge values in heatmaps, and can be further split by whether the heatmap is static or updated during inference. (1) NAR methods with static heatmaps typically take advanced heuristic algorithms with k-opt, such as LKH~\cite{helsgaun2000effective}, as their backbone. For example, NeuroLKH~\cite{xin2021neurolkh} replaces LKH's handcrafted edge-preference prediction rule with a learned one to determine edge candidate sets and search priorities. (2) NAR methods with dynamic heatmaps often utilize MCTS to iteratively update the heatmaps for guiding the k-opt search. In particular, Att-GCN~\cite{fu2021generalize} merges multiple heatmaps from small-scale subgraphs to generate the heatmap for a large-scale instance. DIMES~\cite{qiu2022dimes} incorporates an extra meta-learning-based fine-tuning stage to improve performance. DIFUSCO~\cite{sun2023difusco} introduces a graph-based diffusion framework for modeling the explicit node or edge selection, while UTSP~\cite{min2024unsupervised} eliminates the need for costly labeled datasets via unsupervised learning. Nevertheless, SoftDist~\cite{xia2024position} critically re-evaluates the heatmap-MCTS-k-opt paradigm, particularly questioning the effectiveness of DL-based heatmap generation. This finding highlights fundamental limitations of the current paradigm, underscoring the need for more principled studies.

\begin{table*}[htbp]
  \centering
  \caption{Representative Improvement-based Single-solution-based Large Neighborhood NRSs}
  \tabcolsep=0.07cm
  % \vspace{-0.2cm}
  \begin{threeparttable}
    \begin{tabular}{llllllll}
    \toprule[0.5mm]
    \multicolumn{1}{p{4.5em}}{\textbf{Quaternary}} & \textbf{Generation} & \textbf{Solvable} & \multirow{2}[1]{*}{\textbf{Backbone}} & \textbf{Learning} & \multirow{2}[1]{*}{\textbf{Method}} & \multirow{2}[1]{*}{\textbf{Year}} & \multirow{2}[1]{*}{\textbf{Remarks}} \\
    \textbf{Category} &    \multicolumn{1}{l}{\textbf{Paradigm}}   & \textbf{VRPs} &       & \textbf{Paradigm} &       &       &  \\
    \midrule
          & \multirow{7}[7]{*}{AR} & CVRP, & \multicolumn{1}{p{5.125em}}{Transformer} & RL    & NLNS~\cite{hottung2020neural}  & 2020  & LNS with two handcrafted destroy and \\
       \multirow{9}[3]{*}{Unrestricted}   &       & \multicolumn{1}{p{6.875em}}{SDVRP} &       &       &       &       & \multicolumn{1}{p{22.44em}}{one learned repair criteria.} \\
\cmidrule{3-8}  \multirow{8}[3]{*}{Direct LNS}  &       & CVRP(TW) & \multicolumn{1}{p{5.125em}}{GAT, GRU} & RL    & EGATE~\cite{gao2020learn} & 2020  & \multicolumn{1}{p{22.44em}}{LNS with a learned criteria for both destroy and repair.} \\

\cmidrule{3-8}          &       & CVRP  & Transformer & RL    & \multicolumn{1}{p{4.5em}}{L2I~\cite{lu2019learning}} & 2020  & ILS with both small and large neighborhood search. \\

\cmidrule{3-8}          &       & \multicolumn{1}{p{6.875em}}{TSP, CVRP} & \multicolumn{1}{p{5.125em}}{Transformer} & SL    & L2C-Insert$^\ddagger$~\cite{luo2025learning} & 2025  & \multicolumn{1}{p{25em}}{LNS with a handcrafted destroy and a learned repair criteria.} \\

\cmidrule{3-8}          &       & CVRP(TW), & \multicolumn{1}{p{5.125em}}{Transformer} & RL    & NDS~\cite{hottung2025neural} & 2025  & LNS with a learned destroy and a handcrafted repair criteria. \\
        &       & PCVRP &  &     &  &   &  \\

\cmidrule{2-8}          & NAR   & TSP   & Transformer & SL    & \multicolumn{1}{p{4.5em}}{GenSCO~\cite{li2025generation}} & 2025  & ILS with a generation process for local search. \\
\cmidrule{2-8}          & /   & CVRP(TW)   & GAT, GRU & SL    & \multicolumn{1}{p{4.5em}}{L2Seg~\cite{ouyang2025learning}} & 2026  & A decomposition technique to accelerate iterative solvers. \\
          &    &    & Transformer &    &  &   & \\
    \midrule
       \multirow{13}[7]{*}{Restricted}   & \multirow{7}[6]{*}{AR} & \multicolumn{1}{l}{\multirow{3}[1]{*}{TSP, CVRP}} & \multirow{3}[1]{*}{Transformer} & \multirow{3}[1]{*}{SL}      & \multicolumn{1}{p{4.5em}}{LEHD$^\ddagger$~\cite{luo2023neural}} & 2023  & Random reconstructions of partial solutions via appending. \\
      \multirow{13}[7]{*}{Direct LNS}    &       &       &       &     & \multicolumn{1}{p{4.5em}}{SIL$^\ddagger$~\cite{luo2025boosting}} & 2025  & \multicolumn{1}{p{24em}}{Parallel reconstructions of partial solutions via appending.} \\
          &       &       &       &       & DRHG~\cite{li2025destroy}  & 2025  & LNS with restricted ranges (outside hypernodes). \\
\cmidrule{3-8}          &       & \multicolumn{1}{p{6.875em}}{(PC)TSP, CVRP} & \multicolumn{1}{p{5.125em}}{Transformer} & RL    & LCP~\cite{kim2021learning}   & 2021  & Iterative re-decompositions and revisions. \\
\cmidrule{3-8}     &       & \multicolumn{1}{p{6.875em}}{(A)TSP, OP,} &       &       &       &       &  \\
          &       & \multicolumn{1}{p{6.875em}}{CVRP, OVRP, } & \multicolumn{1}{p{5.125em}}{AGNN,} & RL    & UDC~\cite{zheng2024udc}   & 2024  & Considering the negative impact of sub-optimal  \\
          &       & \multicolumn{1}{p{6.875em}}{(S)PCTSP,} & \multicolumn{1}{p{5.125em}}{Transformer} &       &       &       & dividing policies. \\
          &       & min-max mTSP &       &       &       &       &  \\
\cmidrule{2-8}          & \multirow{4}[6]{*}{/} &
\multirow{1}[5]{*}{CVRP(TW),} &
\multicolumn{1}{p{5.125em}}{Transformer} & RL & L2D~\cite{li2021learning} & 2021 &
Iterative subproblem selection and optimization. \\
\cmidrule{4-8}          &      &    \multirow{1}[1]{*}{VRPMPD}    &
\multicolumn{1}{p{5.125em}}{LSTM,} & RL & RBG~\cite{zong2022rbg} & 2022 &
Iterative re-partitioning, merging, and re-solving. \\
         &       &       &
\multicolumn{1}{p{5.125em}}{Transformer} &   &   &   &  \\

\cmidrule{3-8}          &       & \multicolumn{1}{p{6.875em}}{(A)TSP, CVRP,} & \multicolumn{1}{p{5.125em}}{GNN,} & RL    & GLOP~\cite{ye2024glop}  & 2024  & \multicolumn{1}{p{22.44em}}{An NRS with both NAR and AR paradigms.} \\
          &       & \multicolumn{1}{p{6.875em}}{PCTSP} & \multicolumn{1}{p{5.125em}}{Transformer} &       &       &       &    \\
    \midrule
    \multirow{2}[1]{*}{Indirect LNS} & \multicolumn{1}{l}{\multirow{2}[1]{*}{NAR}} & \multirow{2}[1]{*}{TSP} & \multirow{2}[1]{*}{AGNN} & \multirow{2}[1]{*}{SL} & T2T~\cite{li2024distribution}   & 2023  & Integrating local search in diffuse-and-denoise. \\
          &       &       &       &       & Fast T2T~\cite{li2024fast} & 2024  & Mapping from different noise levels to the optima. \\
    \bottomrule[0.5mm]
    \end{tabular}%
    \begin{tablenotes}
      \footnotesize
      \item[$\ddagger$] The NRS supports multiple inference strategies and currently employs an improvement-based one. For LEHD, SIL, and L2C-Insert, RRC, PRC, and insertion-based local reconstruction are adopted, respectively.
    \end{tablenotes}
  \end{threeparttable}
    % \vspace{-0.6cm}
  \label{tab:LNS}%
\end{table*}%

\subsection{Large Neighborhood Methods}
\label{Section 4.2: LNS}

Large neighborhood methods are grounded in the LNS heuristics~\cite{shaw1998using}, which explore broader solution regions to escape local optima while maintaining manageable computational complexity~\cite{pisinger2019large}. Corresponding NRSs learn different rules to either enhance classical LNS components, such as destroy and repair criteria for perturbation, or to automate the criterion selection. Beyond refining classical LNS, NRSs also introduce novel paradigms such as search in auxiliary latent spaces. These approaches are identified as direct LNS when searching directly on the original solution representation, and indirect LNS when conducted in an auxiliary space.

\subsubsection{Direct LNS}

Direct LNS methods search directly on the original decision space. They can be further categorized by the flexibility of allowed modifications: (a) unrestricted direct LNS permits modifications anywhere in the solution sequence, whereas (b) restricted direct LNS limits modifications to certain predefined positions of the solution.

\paragraph{Unrestricted Direct LNS}

Unrestricted Direct LNS methods are typically built upon the classic destroy-and-repair paradigm of the LNS heuristic, where neighborhoods are implicitly defined by the destroy and repair criteria. In each iteration, the destroy step removes multiple nodes from the complete solution, and the repair step reinserts them sequentially back into the solution for potential improvement. This approach offers two key advantages:  (1) computational scalability~\cite{gendreau2010handbook}, because the number of nodes removed and reinserted (\ie, the perturbation strength) is independent of instance size; and (2) solution quality~\cite{mara2022survey}, as even a small set of nodes, when destroyed and repaired under effective criteria, can lead to promising improvement.

Current related NRSs typically focus on learning effective destroy or repair criteria. For example, NLNS~\cite{hottung2020neural} incorporates two handcrafted destroy criteria and a learned repair criterion, while EGATE~\cite{gao2020learn} employs a single learned rule to both select nodes for removal and determine reinsertion sequences. When applied iteratively, L2C-Insert ~\cite{luo2025learning} can also be regarded as an LNS variant when its learned insertion rule is treated as the repair criterion, complemented by a handcrafted destroy step. The \underline{I}terated \underline{L}ocal \underline{S}earch (ILS) heuristic~\cite{stutzle1999local} further extends LNS by interleaving between large neighborhood perturbation to escape the local region and fine-grained local search to refine the solution. In particular, L2I~\cite{lu2019learning} integrates DL into ILS to select both local search operators and destroy or repair criteria. GenSCO~\cite{li2025generation} perturbs solutions via successive 2-opt moves, as commonly used in heuristics~\cite{penna2013iterated}, and then refines them using a rectified flow model.

The effectiveness of LNS heuristics relies not only on well-designed destroy and repair criteria, but also on rules for controlling the perturbation strength, adapting the criteria, determining the insertion orders, and designing more complex acceptance criteria~\cite{lourencco2019iterated,brandao2020memory,pisinger2018large, gendreau2010handbook}. However, the current NRSs have focused predominantly on learning destroy and repair criteria, leaving other critical rules still largely handcrafted. Therefore, a key future research direction is to automate the design of these rules and to thoroughly investigate their interactions. This holistic design principle is crucial for advancing both this subcategory and NRSs more broadly.

\paragraph{Restricted Direct LNS}
\label{Section 4.2.1.2: Direct LNS (Restricted)}

After the destroy step, typical LNS heuristics encounter scenarios involving partial solutions and unvisited nodes, identical to those faced in construction-based methods. Recent studies have therefore drawn inspiration from single-stage and two-stage construction-based methods to develop new iterative approaches adhering to LNS principles. Some of them iteratively reconstruct partial solutions with single-stage strategies. In contrast, others adopt iterative versions of two-stage methods,  which repeatedly partition the problem and solve subproblems with existing NRSs or heuristics. Though not explicitly framed in classical heuristics, these paradigms can be regarded as position-restricted destroy-and-repair and thus a subcategory of LNS methods.

For the extensions of single-stage methods, the appending LEHD~\cite{luo2023neural} can use a flexible \underline{R}andom \underline{R}e-\underline{C}onstruct (RRC) approach to refine a sampled partial solution at each iteration. From the perspective of LNS, RRC, and its parallel version \underline{P}arallel local \underline{R}e-\underline{C}onstruction (PRC), destroy random node sequences and adopt the learned appending rule as the repair criterion. Integrating this iterative improvement approach into SL training can further reduce the reliance on high-quality solutions~\cite{luo2023neural,luo2025boosting,pirnay2024self} and even enable direct training on large-scale instances~\cite{luo2025boosting}. In addition, DRHG~\cite{li2025destroy} treats partial solutions as hypernodes, which are sequentially appended with unconnected nodes during the repair process.

Instead of designing more powerful subsolvers, the extensions of two-stage methods, such as the so-called ``hierarchical search''~\cite{kim2021learning,zong2022rbg}, ``divide-and-conquer approach''~\cite{cheng2023select,ye2024glop,zheng2024udc}, and ``learning-augmented local search''~\cite{li2021learning}, focus on developing appropriate strategies that leverage existing NRSs or heuristics to achieve better overall performance. For example, LCP~\cite{kim2021learning} employs a seeder policy to generate candidate solutions, which are then optimized in parallel by a reviser that iteratively decomposes and reconstructs them. RBG~\cite{zong2022rbg} decomposes a complete solution into non-overlapping regions, each containing several routes. This division is iteratively updated by a learned rewriter that selects regions to split or merge, after which a generator then generates the routes for the updated regions. Both GLOP~\cite{ye2024glop} and UDC~\cite{zheng2024udc} initially generate a heatmap for decomposition. GLOP partitions the original problem once into sub-TSPs, in which divide-and-conquer steps are further applied. In contrast, UDC employs iterative subproblem re-divisions, and the subproblems are not limited to TSP.

These methods have gained popularity owing to their ability to extend existing construction-based NRSs through iterative refinement. However, like two-stage construction-based methods, the subsolvers generally lack global information, which potentially leads to premature convergence. Moreover, they are often presented merely as extensions of construction-based NRSs, without explicitly acknowledging their LNS nature, leading to insufficient attention to holistic algorithm design. A systematic analysis of LNS heuristics could inspire more principled designs. Promising future research directions include dynamically controlling subproblem sizes, similar to adaptive perturbation degree control in LNS, to balance exploration and exploitation, and identifying suboptimal partial solutions for further improvement while preserving promising ones.

\subsubsection{Indirect LNS}
\label{Section 4.2.3: Indirect LNS}

LNS methods can be generalized to operate in an auxiliary space rather than the original decision space. The auxiliary space is often continuous, enabling gradient-based methods to guide the search. Moreover, operations performed in this space can simultaneously modify multiple parts of a solution, bypassing the sequential node-by-node selection and positioning required in the original decision space. It enables more extensive solution adjustment compared to classical destroy-and-repair perturbations, which typically modify only a small number of edges.

A typical example is to use diffusion models for solving TSPs in an NAR manner~\cite{li2024distribution, li2024fast}. During inference, the forward noising process gradually increases the confidence of extra edges, which turns a feasible solution into an infeasible one with more edges. Conversely, the reverse denoising process decreases the confidence of redundant edges to recover a feasible solution. The stochastic nature of diffusion allows applying an iterative noising-denoising process to produce diverse solutions. Incorporating effective guidance, such as the gradient feedback in T2T for denosing~\cite{li2024distribution}, can further improve the solution quality. A subsequent work, Fast T2T~\cite{li2024fast}, further accelerates denoising via consistency modeling.

Although indirect LNS methods differ from classical LNS, heuristics can still offer valuable insights. Current implementations typically employ a fixed noise schedule during inference. However, as stated earlier, adaptive perturbation strength is crucial for balancing fine-grained search and escape from local optima~\cite{gendreau2010handbook}. Therefore, adaptively adjusting re-noising levels based on search progress could be helpful. Additionally, the greedy decoding is often suboptimal, and more well-designed inference strategies deserve greater attention as in other NAR NRSs. Finally, heatmap-guided search is not the only possible paradigm for indirect LNS. Further work is expected to explore alternative auxiliary search spaces.

\section{Population-Based Methods\\For Improvement}

\label{Section 5: Population-based}

\begin{table*}[ht]
  \centering
  \caption{Representative Improvement-based Population-based NRSs}
  \tabcolsep=0.142cm
  % \vspace{-0.2cm}
    \begin{tabular}{llp{6.19em}lp{4.25em}lll}
    \toprule[0.5mm]
    \textbf{Search} & \textbf{Generation} & \multicolumn{1}{l}{\textbf{Solvable}} & \multirow{2}[1]{*}{\textbf{Backbone}} & \multicolumn{1}{l}{\textbf{Learning}} & \multirow{2}[1]{*}{\textbf{Method}} & \multirow{2}[1]{*}{\textbf{Year}} & \multirow{2}[1]{*}{\textbf{Remarks}} \\
    \textbf{Space} & \textbf{Paradigm} & \multicolumn{1}{l}{\textbf{VRPs}} &       & \multicolumn{1}{l}{\textbf{Paradigm}} &       &       &  \\
    \midrule
    \multirow{2}[1]{*}{Continuous} & \multirow{2}[1]{*}{NAR} & \multirow{2}[1]{*}{TSP, CVRP} & \multicolumn{1}{p{4.565em}}{GRU} & \multicolumn{1}{l}{SL} & CVAE-Opt~\cite{hottung2021learning} & 2021  & Latent space search with DE. \\
     &       & \multicolumn{1}{l}{} & \multicolumn{1}{l}{Transformer} & \multicolumn{1}{l}{RL} & COMPASS~\cite{chalumeau2023combinatorial} & 2023  & Latent space search with CMA-ES. \\
    \midrule
    \multirow{2}[1]{*}{Discrete} & \multirow{2}[1]{*}{NAR} & \multicolumn{1}{p{6.8em}}{(PC)TSP, (S)OP,} & \multicolumn{1}{l}{\multirow{2}[1]{*}{GNN}} & RL & DeepACO~\cite{ye2024deepaco} & 2023  & Learning heuristic measures in ACO with RL. \\
     &       & CVRP(TW) &       & GFlowNet    & GFACS~\cite{kim2024ant} & 2025  & Learning heuristic measures in ACO with GFlowNet. \\
    \bottomrule[0.5mm]
    \end{tabular}%
    % \vspace{-0.6cm}
  \label{tab:population}%
\end{table*}%

Population-based methods maintain and evolve a set of candidate solutions, leveraging collective information from the entire set to guide search~\cite{vidal2014unified,kennedy2006swarm,back1996evolutionary}. In NRSs, these methods can be implemented either by operating directly on the discrete solution space of the original problem, or by transforming solutions into a continuous latent space for optimization, as illustrated in Table~\ref{tab:population}.

For methods that work in the discrete solution space, DeepACO~\cite{ye2024deepaco} and GFACS~\cite{kim2024ant} enhance the classic \underline{A}nt \underline{C}olony \underline{O}ptimization (ACO) by replacing handcrafted heuristic measures for edges (\eg, inverting the length) with learned scoring rules. Unlike various NAR construction-based NRSs confined to TSP and \underline{M}aximum \underline{I}ndependent \underline{S}et (MIS), these approaches inherit the flexibility of meta-heuristics, which can tackle a broader range of COPs.

For methods that work in continuous latent space, CVAE-Opt~\cite{hottung2021learning} utilizes a \underline{V}ariational \underline{A}uto\underline{e}ncoder (VAE) model to learn the distribution of high-quality solutions, then evolves a population in the latent space via \underline{d}ifferential \underline{e}volution (DE)~\cite{price2006differential}. Besides, COMPASS~\cite{chalumeau2023combinatorial} parametrizes a continuous policy distribution and applies \underline{C}ovariance \underline{M}atrix \underline{A}daptation \underline{E}volution \underline{S}trategy (CMA-ES)~\cite{hansen2016cma} to search.

Like classic domain-agnostic meta-heuristics, population-based methods exhibit inherent robustness for problems with complex search spaces. A promising future direction is to adapt them to problems with dynamic environments, where traditional population-based heuristics have demonstrated strong suitability~\cite{gendreau2010handbook}. In addition, given the successful DL-based enhancement of the single-solution-based LKH~\cite{helsgaun2000effective} (as discussed in Section \ref{Section 4.1.2: Sequential}), powerful population-based algorithms like HGS~\cite{vidal2014unified} could likewise be integrated with DL techniques to develop more competitive NRSs.

\section{Experimental Studies}
\label{Section 6: Experiments}

This section investigates the in-problem performance of representative NRSs, with a focus on their zero-shot generalization ability, a topic of significant interest in recent years. The conventional evaluation pipeline is first applied, which emphasizes scalability on synthetic instances and yields promising results. Nevertheless, this pipeline suffers from notable limitations, including a narrow range of test distributions, conflated in- and out-of-distribution comparisons, and inconsistent inference settings. Therefore, a generalization-focused evaluation pipeline is introduced for single-model performance across diverse benchmark instances, with unified inference and complementary metrics. Experimental results under this new pipeline reveal that NRSs trained on narrowly distributed data may be outperformed by even simple construction heuristics such as nearest neighbor and random insertion. This contrast suggests that the conventional pipeline can systematically lead to overly optimistic conclusions. Building on these findings, the advantages of the proposed pipeline are discussed, and principles for method selection are outlined. In particular, learning is argued to remain crucial for NRSs, even when their performance falls short of prior expectations. The implementation details of the experimental studies are available in \url{https://github.com/CIAM-Group/NRS_Survey}.

\subsection{Selected Methods for Comparative Evaluation}

The comparative evaluation incorporates two groups of methods: classical and SOTA heuristics that serve as baselines, and representative NRSs. The selected heuristics, chosen for their efficiency or effectiveness, are briefly introduced below.
\begin{itemize}
    \item \textbf{Nearest Neighbor } A classic construction-based heuristic. At each step, the nearest node to the last node of the partial solution is selected for appending.
    \item \textbf{Random Insertion } A classic construction-based heuristic. At each step, a randomly selected node is inserted at the position that minimizes the increase in cost.
    \item \textbf{LKH-3}~\cite{helsgaun2017extension} ~A single-solution-based SOTA heuristic for TSP, widely adopted as a baseline in prior works.
    \item \textbf{HGS}~\cite{vidal2022hybrid} ~A population-based SOTA heuristic for CVRP, widely adopted as a baseline in prior works.
    \item \textbf{AILS-II}~\cite{maximo2024ails} ~A single-solution-based SOTA heuristic for CVRP, rarely adopted as a baseline in prior works.
\end{itemize}

The selected NRSs comprehensively cover all categories in the proposed taxonomy and are listed in Table~\ref{tab:methods_in_exp}. All inference experiments of NRSs are uniformly conducted on a single NVIDIA GeForce RTX 3090 GPU with 24GB of memory. Specifically, 20 cores of the Intel(R) Xeon(R) Gold 6348 CPU @ 2.60GHz and 40 GB of memory are allocated to each NAR NRS (GFACS, GenSCO, and Fast T2T) for potential calculations on the CPU.

\begin{table}[t]
  \centering
  \caption{Selected NRSs for Comparative Evaluation}
  % \vspace{-0.2cm}
  \begin{threeparttable}
  \tabcolsep=0.11cm
    \begin{tabular}{ccccl}
    \toprule[0.5mm]
    \multicolumn{4}{c}{Category}  & \multicolumn{1}{c}{\multirow{2}[4]{*}{Method}} \\
    \cmidrule{1-4}
    Primary & Secondary & Tertiary & Quaternary &  \\
    \midrule
    \multirow{11}[6]{*}{\rotatebox[origin=c]{90}{Construction}} & \multirow{10}[5]{*}{Single-stage} & \multirow{9}[1]{*}{Appending} & \multirow{9}[1]{*}{/} & BQ~\cite{drakulic2024bq} \\
          &       &       &       & \multicolumn{1}{p{4.565em}}{LEHD$^{\dagger}$~\cite{luo2023neural}} \\
          &       &       &       & \multicolumn{1}{p{4.565em}}{SIL$^{\dagger}$~\cite{luo2025boosting}} \\
          &       &       &       & ICAM~\cite{zhou2024instance} \\
          &       &       &       & ELG~\cite{gao2024towards} \\
          &       &       &       & INViT~\cite{fang2024invit} \\
          &       &       &       & L2R~\cite{zhou2025l2r} \\
          &       &       &       & DGL~\cite{xiao2025dgl} \\
          &       &       &       & ReLD~\cite{huang2025rethinking} \\
    \cmidrule{3-5}
          &       & Insertion & /     & \multicolumn{1}{p{7.2em}}{L2C-Insert$^{\dagger}$~\cite{luo2025learning}} \\
    \cmidrule{2-5}
          & Two-stage & /     & /     & H-TSP~\cite{pan2023htsp} \\
    \midrule
    \multirow{9}[14]{*}{\rotatebox[origin=c]{90}{Improvement}} & \multirow{8}[12]{*}{Single-solution} & Small & Immediate & DACT~\cite{ma2021learning} \\
    \cmidrule{4-5}
          &       & Neighborhood & Sequential & NeuOpt~\cite{ma2024learning} \\
    \cmidrule{3-5}
      &       &       
      & Unrestricted
      & \multicolumn{1}{p{7.2em}}{L2C-Insert$^{\ddagger}$~\cite{luo2025learning}} \\
          &       &       & Direct LNS & GenSCO~\cite{li2025generation} \\
    \cmidrule{4-5}
      &       & Large 
      & \multirow{3}{*}{\makecell[c]{Restricted \\ Direct LNS}}
      & \multicolumn{1}{p{4.565em}}{LEHD$^{\ddagger}$~\cite{luo2023neural}} \\
      &       & Neighborhood 
      & 
      & SIL$^{\ddagger}$~\cite{luo2025boosting} \\
      &       &       
      & 
      & DRHG~\cite{li2025destroy} \\
    \cmidrule{4-5}
          &       &       & Indirect LNS & Fast T2T~\cite{li2024fast} \\
    \cmidrule{2-5}
          & Population & /     & /     & GFACS~\cite{kim2024ant} \\
    \bottomrule[0.5mm]
    \end{tabular}
    \begin{tablenotes}
      \footnotesize
      \item[$\dagger$] The NRS supports multiple inference strategies and currently employs a construction-based greedy inference.
      \item[$\ddagger$] The NRS supports multiple inference strategies and currently employs an improvement-based one. For LEHD, SIL, and L2C-Insert, RRC, PRC, and insertion-based local reconstruction are adopted, respectively.
    \end{tablenotes}
  \end{threeparttable}
  % \vspace{-0.6cm}
  \label{tab:methods_in_exp}
\end{table}

\subsection{Experiment on Conventional Evaluation Pipeline}

\begin{table}[t]
  \centering
  \caption{Experimental Results of Conventional Evaluation Pipeline}
  \tabcolsep=0.02cm
  % \vspace{-0.2cm}
  \begin{threeparttable}
    \begin{tabular}{l|cc|cc|cc}
    \toprule[0.5mm]
    \multicolumn{1}{c|}{\multirow{2}[1]{*}{Method}} & \multicolumn{2}{c|}{TSP 100} & \multicolumn{2}{c|}{TSP 1K} & \multicolumn{2}{c}{TSP 10K} \\
    % \cline{2-7}
     & Gap & Time & Gap & Time & Gap & Time \\
    \midrule
    LKH-3 & 0.000\% & 10.97m & 0.000\% & 5.69m & 0.000\% & 49.34m \\
    \midrule
    Nearest Neighbor & 24.722\% & 6.72s & 25.022\% & 0.97s & 23.864\% & 2.29s \\
    Random Insertion & 9.672\% & 2.04s & 13.096\% & 0.46s & 13.966\% & 4.64s \\
    \midrule
    $^\uparrow$BQ greedy   & 0.348\% & 1.13m & 2.294\% & 1.19m & /     & / \\
    $^\uparrow$LEHD$^*$ greedy & 0.576\% & 26.84s & 3.116\% & 1.64m & /     & / \\
    $^\parallel$SIL$^*$ greedy & /     & /     & 1.952\% & 29.12s &  { {4.061\%}} & 6.06m \\
    $^\uparrow$ICAM aug $\times 8$ &  { {0.147\%}} & 44.66s &  { {1.647\%}} & 3.93m & /     & / \\
    $^\uparrow$ELG aug $\times 8$ & 0.224\% & 3.02m & /     & /     & /     & / \\
    $^\uparrow$INViT-3V aug$^\dagger$ & 1.419\% & 32.44m & 5.154\% & 5.52m & 6.678\% & 1.27h \\
    $^\uparrow$L2R greedy  & /     & /     & 4.494\% &  { {6.48s}} & 4.824\% &  { {1.07m}} \\
    $^\uparrow$DGL aug$^\dagger$ & 0.609\% & 14.16m & 2.714\% & 1.42m & 6.792\% & 10.61m \\
    $^\uparrow$L2C-Insert$^*$ greedy & 0.458\% & 1.24m & 4.756\% & 32.98s & 7.760\% & 1.11m \\
    \midrule
    $^\parallel$H-TSP & /     & /     & 6.673\% & 46.59s & 8.329\% & 50.92s \\
    \midrule
    $^\parallel$DACT T=10K & 0.379\% & 2.05h & /     & /     & /     & / \\
    $^\parallel$NeuOpt T=10K & 0.018\% & 1.47h & /     & /     & /     & / \\
    \midrule
    $^\uparrow$L2C-Insert$^*$ T=1K &  { {0.0001\%}} & 12.01h & 0.485\% & 1.21h & 2.086\% & 15.86m \\
    $^\parallel$GenSCO 2-opt & 0.0003\% &  { {1.96m}} &  { {0.033\%}} &  { {6.76m}} & /     & / \\
    $^\uparrow$LEHD$^*$ RRC1K & 0.002\% & 2.36h & 0.729\% & 7.49h & /     & / \\
    $^\parallel$SIL$^*$ PRC1K & /     & /     & 0.375\% & 3.47h & 1.824\% & 5.19h \\
    $^\uparrow$DRHG T=1K & 0.0003\% & 7.12h & 0.420\% & 3.89h &  { {1.802\%}} & 1.05h \\
    $^\parallel$Fast T2T$^\ddagger$ T\textsubscript{s}=5, T\textsubscript{g}=5 \quad & 0.030\% & 37.29m & 0.589\% & 9.03m & /     & / \\
    \midrule
    $^\parallel$GFACS T=10, K=100 & /     & /     & 2.615\% & 3.14h & /     & / \\
    \bottomrule[0.5mm]
    \end{tabular}%
    \begin{tablenotes}
      \footnotesize
      \item[$\uparrow$] The marked NRS reports generalization performance of single models, typically trained on the smallest-scale instances reported and tested on larger-scale, all with uniform node distribution.
      \item[$\parallel$] The marked NRS reports in-distribution performance of multiple models on specific scales with uniform node distribution.
      \item[$*$] The NRS supports multiple inference strategies.
      \item[$\dagger$] For INViT and DGL,  data augmentation factors vary by scale. Besides, the INViT model with kNN size of (65, 50, 35) is adopted instead of the unavailable (50, 35, 15) reported in the original study.
      % \item[] For GenSCO, $C=10$ is adopted in TSP 100 and $C=160$ is adopted in TSP 1K, according to the original paper.
      \item[$\ddagger$] For Fast T2T, the additional 2-opt improvement as a post process is adopted only in TSP 1K, according to the original paper.
    \end{tablenotes}
  \end{threeparttable}
  % \vspace{-0.6cm}
  \label{tab:exp_conventional}%
\end{table}%

\subsubsection{Experimental Purpose}

This pipeline generally evaluates NRSs on synthetic instances with specific scales, node distributions, and optional constraint tightness~\cite{luo2025rethinking}. Among these aspects,  scalability is the most widely studied one and is also the primary focus of this experiment. It is important to note, however, that scalability is not equivalent to generalization, which will be discussed in detail in Section~\ref{Section 6.4.1 Pipeline Advances}.

\subsubsection{Experimental Settings}

\paragraph{Problem and Instance Setting}

TSP unrelated to constraint tightness is considered due to the lack of a unified setting in the literature. For scale and node distribution, the evaluation follows common practice by testing on uniformly distributed instances at scales of 100, 1K, and 10K. All instances are drawn from the generated datasets of SIL~\cite{luo2025boosting}.

\paragraph{Metrics and Inference}

Two metrics are reported for each method: the optimality gap (Gap) and the total inference time (Time). Specifically, the optimality gap measures the discrepancy between the obtained solutions and the best-known solutions, provided by the LKH-3 heuristic, as is common practice in AM~\cite{kool2018attention}. For NRSs, the released implementations and pretrained models are adopted. Note that each NRS is evaluated only under a specific configuration on instances with corresponding sizes reported in the original studies. Results for unreported conditions are denoted by ``/''.

\subsubsection{Performance Evaluation}

According to Table~\ref{tab:exp_conventional}, \textbf{NRSs exhibit promising performance under the conventional pipeline}. \textbf{For construction-based NRSs, all of them outperform simple heuristics (nearest neighbor and random insertion) within their respective categories.} Specifically, ICAM achieves strong in- and out-of-distribution results. Besides, L2R maintains competitive performance while reducing inference time by approximately an order of magnitude. For large-scale instances with 10K nodes, where only a few construction-based methods are evaluated, SIL (Greedy) delivers the best performance. In contrast, H-TSP generally underperforms single-stage counterparts, falling short of the expected two-stage advantages on larger instances. \textbf{For improvement-based NRSs, most of them achieve competitive performance close to that of the advanced heuristic LKH-3.} For example, GenSCO with 2-opt achieves strong in-distribution results within a short runtime. In addition, among the limited NRSs tested at 10K, DRHG performs best, achieving slightly better performance than SIL (PRC) while using only about one-fifth of its inference time. Nevertheless, GFACS is outperformed by several construction-based methods (BQ, SIL (Greedy), and ICAM) at the scale of 1K.

\subsection{Experiment on the Proposed Evaluation Pipeline}

\subsubsection{Experimental Purpose}

The conventional evaluation pipeline has several limitations. First, its testing distributions are limited in scope, typically restricted to specific scales and node distributions~\cite{zhou2023towards}. This restricted coverage poorly represents real-world scenarios. Moreover, the parameterized synthetic instance generators can bias performance toward certain training distributions. Second, it does not distinguish the evaluation of single-model generalization performance (\ie, one model applied to all test instances) and multi-model in-distribution performance (\ie, separate models trained and tested per problem scale). Finally, the inference settings are typically inconsistent across different methods. In short, this pipeline inherently favors NRSs whose DL models overfit to the training distribution and report multi-model in-distribution performance, therefore introducing systematic evaluation bias.

To address these issues, a new evaluation pipeline is introduced. It centers on the zero-shot in-problem generalization, which has been the primary focus of advanced NRSs in recent years and therefore serves as a representative indicator of progress in the field. Under this pipeline, NRSs are benchmarked on diverse instances that more faithfully reflect the irregular conditions of real-world production and logistics. Distributional biases inherent in synthetically generated instances, particularly those with uniformly distributed nodes, are avoided. In addition, the inference settings are consistently standardized across all evaluated NRSs.

\begin{table*}[htbp]
  \centering
  \caption{Experimental Results of the Proposed Evaluation Pipeline}
  % \vspace{-0.2cm}
  \tabcolsep=0.185cm
  \begin{threeparttable}
    \begin{tabular}{c|l|ccc|ccc|ccc|cc}
    \toprule[0.5mm]
          & \multicolumn{1}{c|}{\multirow{2}[2]{*}{Method}} & \multicolumn{3}{c|}{(0,1K)} & \multicolumn{3}{c|}{[1K, 10K)} & \multicolumn{3}{c|}{[10K, 100K]} & \multicolumn{2}{c}{Total} \\
\cmidrule{3-13}          &       & Gap   & Time  & \multicolumn{1}{l|}{Solved} & Gap   & Time  & \multicolumn{1}{l|}{Solved} & Gap   & Time  & \multicolumn{1}{l|}{Solved} & Gap   & Solved \\
\cmidrule{1-13}    \multirow{24}[12]{*}{\rotatebox[origin=c]{90}{TSP}} & Nearest Neighbor & 25.29\% & 0.01s & 69/69 & 26.66\% & 0.29s & 109/109 & 25.01\% & 22.60s & 50/50 & 25.88\% &  {{228}/228} \\
          & Random Insertion & 10.60\% &  {{0.00s}} & 69/69 & 15.32\% &  {{0.05s}} & 109/109 & 16.37\% &  {{8.93s}} & 50/50 & 14.12\% &  {{228}/228} \\
          & LKH-3\phantom{$^{\downarrow}$} t=n/3, runs=1 & 0.00\% & 7.88s & 69/69 & 0.01\% & 631.34s & 109/109 & 0.08\% & 14600.50s & 50/50 & 0.03\% &  {{228}/228} \\
          & LKH-3$^{\downarrow}$ t=n/3, runs=1 & 0.00\% & 9.25s & 69/69 & 0.01\% & 600.36s & 109/109 & 0.05\% & 10800.24s & 50/50 &  {{0.02\%}} &  {{228}/228} \\
\cmidrule{2-13}          & BQ    & 5.00\% & 2.51s & 68/69 & 19.03\% & 22.74s & 92/109 & 52.00\% & 187.81s & 4/50  & 14.02\% & 164/228 \\
          & LEHD$^*$ greedy & 4.85\% & 1.01s & 69/69 & 20.13\% & 68.27s & 106/109 & 49.35\% & 1386.01s & 11/50 & 16.19\% & 186/228 \\
          & SIL$^*$ greedy & 8.64\% & 1.69s & 69/69 & 9.83\% & 17.69s & 109/109 & 11.11\% & 430.73s & 50/50 & 9.75\% &  {{228}/228} \\
          & ICAM  & 6.53\% & 0.25s & 69/69 & 16.62\% & 21.57s & 109/109 & 21.34\% & 1050.33s & 19/50 & 13.54\% & 197/228 \\
          & ELG   & 6.05\% & 0.63s & 69/69 & 18.14\% & 88.12s & 108/109 & 21.65\% & 940.35s & 6/50  & 13.70\% & 183/228 \\
          & INViT-3V   & 7.93\% & 2.77s & 69/69 & 12.08\% & 49.03s & 109/109 & 11.52\% & 1079.50s & 42/50 & 10.67\% & 220/228 \\
          & L2R   & 5.89\% & 1.60s & 69/69 & 9.22\% & 15.55s & 109/109 & 8.52\% & 153.11s & 50/50 &  {{8.06\%}} &  {{228}/228}   \\
          & DGL   & 6.53\% & 1.17s & 69/69 & 11.32\% & 11.67s & 109/109 & 11.14\% & 58.62s & 25/50 & 9.67\% & 203/228 \\
          & L2C-Insert$^*$ greedy & 4.39\% & 1.51s & 69/69 & 18.12\% & 15.34s & 109/109 & 30.94\% & 145.17s & 50/50 & 16.77\% &  {{228}/228}   \\
\cmidrule{2-13}          & H-TSP & 6.16\% &  {{0.61s}} & 36/69 & 11.62\% &  {{3.15s}} & 100/109 & 12.29\% &  {{21.44s}} & 40/50 & 10.65\% & 176/228 \\
\cmidrule{2-13}          & DACT T=1K & 16.37\% & 39.84s & 69/69 & 26.58\% & 261.73s & 83/109 & /     & /     & 0/50  & 21.94\% & 152/228 \\
          & NeuOpt T=1K & 19.90\% & 81.22s & 46/69 & /     & /     & 0/109 & /     & /     & 0/50  & 19.90\% & 46/228 \\
\cmidrule{2-13}          & L2C-Insert$^*$ T=1K & 1.08\% & 381.55s & 69/69 & 9.80\% & 479.98s & 109/109 & 29.15\% & 615.19s & 50/50 & 11.41\% &  {{228}/228}   \\
          & GenSCO & 14.56\% & 23.19s & 68/69 & 35.46\% & 677.31s & 104/109 & 35.17\% & 14304.70s & 25/50 & 28.21\% & 197/228 \\
          & LEHD$^*$ RRC1K & 1.73\% & 498.40s & 69/69 & 10.87\% & 1634.51s & 109/109 & 24.02\% & 2769.86s & 9/50  & 8.13\% & 187/228 \\
          & SIL$^*$ PRC1K & 0.80\% & 883.87s & 69/69 & 2.58\% & 2880.46s & 109/109 & 4.55\% & 4933.45s & 50/50 & 2.47\% &  {{228}/228} \\
          & DRHG T=1K & 0.10\% & 769.53s & 69/69 & 1.46\% & 2857.55s & 109/109 & 4.46\% & 3004.28s & 50/50 &  {{1.71\%}} &  {{228}/228}   \\
          & Fast T2T T\textsubscript{s}=10, T\textsubscript{g}=10 & 10.46\% & 1.41s & 45/69 & /     & /     & 0/109 & /     & /     & 0/50  & 10.46\% & 45/228 \\
\cmidrule{2-13}          & GFACS$^{\dagger}$ T=100, K=100 & 31.64\% & 166.94s & 66/69 & 86.77\% & 2601.13s & 22/109 & /     & /     & 0/50  & 45.42\% & 88/228 \\
          & GFACS$^{\ddagger}$ T=100, K=100 & 0.72\% & 174.21s & 69/69 & 3.76\% & 9142.93s & 83/109 & /     & /     & 0/50  & 2.38\% & 152/228 \\
    \midrule
    \midrule
    \multirow{22}[11]{*}{\rotatebox[origin=c]{90}{CVRP}} & Nearest Neighbor & 21.17\% & 0.03s & 99/99 & 15.18\% & 1.08s & 5/5   & 11.80\% & 14.63s & 6/6   & 20.39\% &  {{110}/110} \\
          & Random Insertion & 75.00\% & 0.00s & 36/99 & /     & /     & 0/5   & /     & /     & 0/6   & 75.00\% & 36/110 \\
          & HGS t=n/3 & 0.29\% & 111.24s & 99/99 & 3.59\% & 1428.41s & 5/5   & 7.86\% & 6926.41s & 6/6   & 0.85\% &  {{110}/110} \\
          & AILS-II t=n/3 & 0.57\% & 133.95s & 99/99 & 1.58\% & 1388.42s & 5/5   & 1.58\% & 5646.48s & 6/6   &  {{0.68\%}} &  {{110}/110} \\
\cmidrule{2-13}          & BQ   & 8.87\% & 3.63s & 99/99 & 20.28\% & 39.97s & 5/5   & 41.52\% & 202.69s & 5/6   & 10.89\% & 109/110 \\
          & LEHD$^*$ greedy & 11.25\% & 1.53s & 98/99 & 19.22\% & 99.43s & 5/5   & 32.80\% & 852.02s & 2/6   & 12.04\% & 105/110 \\
          & SIL$^*$ greedy & 40.04\% & 2.48s & 65/99 & 16.09\% & 26.86s & 5/5   & 10.81\% & 146.83s & 6/6   & 36.16\% & 76/110 \\
          & ICAM  & 5.00\% & 0.42s & 99/99 & 11.69\% & 32.32s & 5/5   & /     & /     & 0/6   & 5.32\% & 104/110 \\
          & ELG   & 8.03\% & 1.29s & 99/99 & 18.51\% & 30.21s & 5/5   & 29.38\% & 133.08s & 2/6   & 8.93\% & 106/110 \\
          & INViT-3V   & 13.15\% & 4.72s & 99/99 & 19.03\% & 77.33s & 5/5   & 23.91\% & 496.25s & 5/6   & 13.91\% & 109/110 \\
          & L2R   & 8.16\% & 2.49s & 99/99 & 11.62\% & 23.99s & 5/5   & 11.08\% & 97.12s & 6/6   & 8.48\% &  {{110}/110}   \\
          & DGL   & 15.27\% & 2.22s & 99/99 & 17.96\% & 22.60s & 5/5   & 18.69\% & 78.64s & 5/6   & 15.55\% & 109/110 \\
          & ReLD  & 4.10\% & 0.41s & 99/99 & 10.22\% & 5.29s & 5/5   & 11.27\% & 28.87s & 3/6   &  {{4.58\%}} & 107/110 \\
          & L2C-Insert$^*$ greedy & 6.87\% & 2.73s & 99/99 & 22.37\% & 616.75s & 5/5   & 49.41\% & 5525.71s & 2/6   & 8.40\% & 106/110 \\
\cmidrule{2-13}          & DACT T=1K & 16.42\% & 246.51s & 74/99 & 17.70\% & 479.82s & 1/5   & /     & /     & 0/6   & 16.44\% & 75/110 \\
          & NeuOpt T=1K & 26.93\% & 571.14s & 36/99 & /     & /     & 0/5   & /     & /     & 0/6   & 26.93\% & 36/110 \\
\cmidrule{2-13}          & L2C-Insert$^*$ T=1K & 3.21\% & 344.05s & 99/99 & 18.87\% & 6166.21s & 5/5   & 44.29\% & 32754.86s & 2/6   & 4.72\% & 106/110 \\
          & LEHD$^*$ RRC1K & 3.58\% & 796.15s & 99/99 & 11.73\% & 2043.74s & 5/5   & 21.98\% & 2820.28s & 2/6   & 4.32\% & 106/110 \\
          & SIL$^*$ PRC1K  & 21.38\% & 1307.97s & 99/99 & 8.28\% & 3471.69s & 5/5   & 7.40\% & 4251.88s & 6/6   & 20.02\% &  {{110}/110}   \\
          & DRHG T=1K & 11.11\% & 1114.60s & 99/99 & 17.95\% & 2529.12s & 5/5   & 16.95\% & 5376.37s & 6/6   & 11.74\% &  {{110}/110}   \\
\cmidrule{2-13}          & GFACS$^{\dagger}$ T=100, K=100 & 36.83\% & 437.38s & 99/99 & 34.09\% & 9654.33s & 3/5   & /     & /     & 0/6   & 36.75\% & 102/110 \\
          & GFACS$^{\ddagger}$ T=100, K=100 & 2.60\% & 405.81s & 99/99 & 7.65\% & 14884.94s & 4/5   & /     & /     & 0/6   & 2.80\% & 103/110 \\
    \bottomrule[0.5mm]
    \end{tabular}%
    \begin{tablenotes}
      \footnotesize
      \item[$\downarrow$] $Initial\_Period$ is set as 1K, rather than the original value DIMENSION/2.
      \item[$*$]  The NRS supports more than one inference strategy.
      \item[$\dagger$] GFACS without local search at the last generation. The output solution is the best of the population at the last generation.
      \item[$\ddagger$] The original version of GFACS. The output solution is the best in history (all after local search). 
      
    \end{tablenotes}
    \end{threeparttable}
    % \vspace{-0.5cm}
  \label{tab:exp_proposed}%
\end{table*}%

\subsubsection{Experimental Settings}

\paragraph{Problem and Instance Setting}

The proposed evaluation pipeline assesses NRSs on representative TSP and CVRP. The test instances are drawn from benchmarks and challenge sets, covering diverse data distributions, with scales in $(0, 100\text{K]}$, and specific edge-weight types (\textbf{EUC\_2D} or \textbf{CEIL\_2D}) to ensure integer Euclidean distance matrices. All selected instances have available \underline{b}est \underline{k}nown \underline{s}olutions (BKS) and do not impose additional constraints, such as fixed route numbers or duration limits. The composition of the test instances is detailed as follows. 

\begin{itemize}

    \item \textbf{TSPLIB}~\cite{reinelt1991tsplib} ~a famous dataset with TSP instances from various sources. 77 \textbf{EUC\_2D} instances and 4 \textbf{CEIL\_2D} are included. Note that the \textbf{EUC\_2D} instance $linhp318$ is excluded due to a fixed-edge constraint.

    \item \textbf{National } a dataset with 27 \textbf{EUC\_2D} TSP instances for countries, based on data from the National Imagery and Mapping Agency. All the instances are included.

    \item \textbf{VLSI } a dataset with 102 \textbf{EUC\_2D} TSP instances of industrial applications of the very large-scale integration design from the Bonn Institute. Note that 4 instances (\textit{SRA104815}, \textit{ARA238025}, \textit{LRA498378}, \textit{LRB744710}) with over 100K nodes are excluded.

    \item \textbf{Dataset of The 8th DIMACS Implementation Challenge (TSP) } a dataset comprises a selection of instances from the TSPLIB library, supplemented by generated instances. To avoid redundant instances and to satisfy distributional diversity and edge-weight-type consistency, only the 22 generated \textbf{EUC\_2D} instances with clustered nodes are included. Note that the instance \textit{C316k.0} with over 100K clustered nodes is excluded.

    \item \textbf{CVRPLIB}~\cite{uchoa2017new} ~a famous dataset with 14 sets of CVRP instances from several academic literature and real-world applications. The library encompasses the adopted open-source instances in the \textbf{12th DIMACS Implementation Challenge (CVRP)}. 100 \textbf{EUC\_2D} instances of Set X~\cite{uchoa2017new} and 10 of Set AGS~\cite{arnold2019efficiently} are included.

\end{itemize}

\paragraph{Metrics}

The solvers are evaluated from the following three perspectives:

\begin{itemize}
    \item \textbf{Effectiveness } a solver's ability to maintain high performance across out-of-distribution instances. It is measured by the average gap relative to the BKSs.

    \item \textbf{Efficiency } a solver's ability to solve the instances in a reasonable time. It is measured by the average computational time a solver requires to output the solutions.

    \item \textbf{Reliability } a solver's ability to successfully solve instances within the current scope. It is measured by the number of instances a solver can handle before failure, where ``failure'' encompasses \underline{O}ut-\underline{o}f-\underline{M}emory (OOM) errors, performance breakdowns (\ie, gaps exceeding 100\%~\cite{Choi2025towards,luo2025boosting}), or timeouts (per-instance runtime beyond 36,000s for NRSs). 
    
\end{itemize}

Results are reported separately for three instance scale groups: small ((0,1\text{K})), medium ([1\text{K},10\text{K})], and large ([10\text{K},100\text{K}]). The overall aggregated results are also provided. Results for unsolvable conditions are denoted by ``/''.

\paragraph{Inference}

For advanced heuristics, the termination criterion follows common practice in the heuristic literature~\cite{maximo2024ails}, where the time budget is set proportional to the instance size. To align with the inference time of NRSs, this multiple is set to one-third. Furthermore, the runtime required to reach the current best solution or the BKS is reported. Results are averaged over 10 independent runs.

For NRSs, to facilitate a direct and fair comparison, all methods are evaluated with greedy inference. In other words, special decoding strategies (\eg, beam search) are deliberately excluded, while data augmentation and fine-tuning techniques (\eg, active search~\cite{bello2016neural,hottung2021efficient}) are deactivated. Unless otherwise specified, additional operator-based local search processes (\eg, 2-opt) are also disabled to preserve experimental fairness and prevent possible shifts in categories of NRSs. All other configurations are kept at their method-specific defaults. For improvement-based methods, the number of iterations is set as the maximum value specified in the original configurations. No additional training is conducted during evaluation. Instead, all publicly available pretrained models (trained on instances with specific scales and uniform node distribution) are tested, and the reported result for each NRS corresponds to the best-performing one, selected by prioritizing reliability first and effectiveness second. 
% Complete results are provided in Tables \textcolor{blue}{XII} and \textcolor{blue}{XIII} in Appendix \textcolor{blue}{A}.
Complete results are provided in Tables~\ref{tab:exp_proposed_TSP_details} and~\ref{tab:exp_proposed_CVRP_details} in Appendix~\ref{Appendix: Exp_Details}.

\subsubsection{Performance Evaluation}

Overall, the results presented in Table~\ref{tab:exp_proposed} lead to conclusions fundamentally different from those under the conventional evaluation pipeline.

\textbf{Overall Performance of NRSs } Under the proposed evaluation pipeline, \textbf{NRSs generally underperform SOTA heuristics in both effectiveness and efficiency}, with the gap widening as the problem size increases. Even with comparable runtime, improvement-based NRSs still fall short of SOTA heuristics across all scales. \textbf{In terms of reliability, only a few NRSs (L2R, SIL (PRC), and DRHG) can successfully solve all TSP and CVRP instances}, among which L2R is the only construction-based method. For the remaining methods, only SIL (Greedy) and L2C-Insert (with both greedy and iterative inference) manage to solve every TSP instance. Notably, all successful cases discussed above benefit from techniques for search space reduction (discussed in Section~\ref{Section 7.1: In-problem}). These results indicate a narrow solvable range of current NRSs.

\textbf{Performance of Construction-based NRSs } \textbf{The performance of construction-based NRSs is less encouraging than that indicated by the conventional evaluation pipeline.} In terms of effectiveness, many single-stage NRSs underperform simple heuristics from the same subcategory. For CVRP, the effectiveness of appending NRS SIL (Greedy) deteriorates at small and medium scales, whereas BQ, LEHD (Greedy), ELG, INViT, and DGL degrade on medium- and large-scale instances. All of these NRSs fall short of the nearest neighbor heuristic. Similarly, for TSP, the insertion NRS L2C-Insert (Greedy) is outperformed by random insertion on medium- and large-scale instances. Nevertheless, a few methods stand out: L2R achieves strong effectiveness and reliability on both problems, and ReLD attains competitive effectiveness on CVRP despite limited reliability on large instances. In terms of efficiency, L2C-Insert (Greedy) runs slower than other single-stage methods on CVRP because its released implementation evaluates all unvisited nodes, rather than restricting attention to the nearest one as described in the paper. As the only two-stage NRS, H-TSP benefits from its architecture to achieve inference times comparable to those of the nearest neighbor heuristic while maintaining stable effectiveness on TSP. Nevertheless, its reliability is not strong even in small-scale instances.

\textbf{Performance of Improvement-based NRSs } The performance of improvement-based NRSs is mixed. \textbf{Among single-solution-based NRSs, LNS methods (especially the direct ones) generally exhibit superior effectiveness and reliability compared to the small neighborhood counterparts}, consistent with their recognized advantages of escaping local optima. For example, DRHG approaches LKH-3's effectiveness on TSP across scales. Nevertheless, a few LNS methods exhibit effectiveness deterioration on large-scale instances (LEHD (RRC) and L2C-Insert (Iteration) for both, DRHG for CVRP), and SIL (PRC) degrades on small-scale CVRP instances. In all these deterioration cases, they perform worse than at least one of the two simple construction-based heuristics.
Besides, L2C-Insert (Iteration) remains inefficient as in its construction-based version. In addition, GenSCO and Fast T2T underperform other large neighborhood NRSs across all evaluation aspects. Fast T2T adopts a distance-based insertion strategy similar to random insertion and achieves comparable effectiveness, suggesting it fails to effectively leverage information from out-of-distribution instances. Besides, GenSCO exhibits a performance drop using greedy decoding without explicitly incorporating distance information. This observation suggests that, for the distance-driven TSP, additional spatial bias during element selection remains important for current NAR NRSs. Lastly, the two small-neighborhood DACT and NeuOpt incur high computational cost in AR node-pair selection and insufficient convergence from per-step sampling, resulting in limited effectiveness and reliability.

For population-based NRSs, two variants of GFACS are evaluated: (1) the original version, which returns the best solution over the entire run, and (2) a variant, which disables local search in the final iteration and outputs the best solution from that iteration. The latter variant aligns with the inference settings of other NRSs, and follows the original study's motivation that local search primarily facilitates convergence during training. The variant shows degraded effectiveness and reliability, indicating that \textbf{edge weights, shaped by learned edge-preference rules and iterative population dynamics, still provide insufficient guidance for constructing high-quality solutions for instances across diverse distributions}.

\subsection{Discussions}

\subsubsection{Advantages of the Proposed Evaluation Pipeline}
\label{Section 6.4.1 Pipeline Advances}

Compared with the conventional evaluation pipeline, the proposed pipeline offers several advantages in the following aspects:

\begin{itemize}
    
    \item \textbf{Purpose of Evaluation } The proposed pipeline is designed specifically to assess the zero-shot generalization performance of NRSs. It enforces a consistent single-model evaluation across methods, thereby enhancing comparability and strengthening the validity of conclusions. In contrast, the conventional pipeline primarily focuses on scalability without clearly specifying whether evaluation is conducted under a single- or multi-model condition. Consequently, this ambiguity makes fair comparisons difficult, as some results reflect generalization from a single model, while others report purely in-distribution performance of multiple models, each evaluated only on its matched distribution.

    \item \textbf{Instance Selection } The proposed pipeline draws test instances from well-known benchmarks and challenge sets, rather than ad-hoc synthetic distributions. It thus enables a more comprehensive evaluation across diverse distributions while mitigating generator-induced distribution shifts that could bias results, preserving fairness and comparability. Notably, although prior work sometimes reports results on TSPLIB or CVRPLIB, the instances are often restricted to selected scale ranges or specific subsets, which may introduce selection bias. In contrast, the proposed pipeline uses a broader instance pool, providing a more robust assessment.

    \item \textbf{Inference and Metrics } For inference, the pipeline enforces a greedy decoding setting and avoids arbitrary add-on enhancements or parameter tuning on specific instances, thereby ensuring a more equitable comparison across methods. For the metrics, in addition to effectiveness, the proposed pipeline explicitly introduces reliability as a complementary metric, enabling a more comprehensive evaluation of algorithmic performance.

\end{itemize}

\subsubsection{Principles for Method Selection}

To ensure fair and informative comparisons, two principles for selecting NRSs and baseline heuristics are followed.

\textbf{In-category Comparison of NRSs } The primary goal of NRS experiments is to demonstrate the effectiveness of the learned heuristics. However, NRSs from different categories may rely on distinct heuristic frameworks, each requiring different levels of domain knowledge and computational resources. Therefore, cross-category comparisons may fail to accurately reflect the specific contribution of a DL model to the overall performance. For this reason, comparisons and conclusions are restricted as much as possible to NRSs and traditional heuristics that belong to the same category.

\textbf{Baseline Heuristic Selection } Results under the proposed pipeline highlighted a performance gap between many NRSs and traditional heuristics, contradicting several existing claims that such NRSs can outperform SOTA heuristic methods~\cite{li2025generation,luo2025boosting,zong2022rbg}. In many of those studies, comparisons are conducted either under settings that disadvantage heuristics or against relatively weak heuristic baselines. For the former cases, when heuristics are allowed a shorter initial period~\cite{fu2023hierarchical} under the same time budget, they can achieve better effectiveness and efficiency on large-scale instances, as shown in Table~\ref{tab:exp_proposed}. For the latter cases, widely-used CVRP baselines in the NRS literature, such as HGS and LKH-3, are not selected for comparison in recent heuristics literature~\cite{maximo2024ails,maximo2021hybrid,christiaens2020slack}. Advanced heuristic methods, such as AILS-II, can achieve more competitive performance on medium- and large-scale instances where NRSs are often claimed to outperform traditional heuristics. The evidence above indicates that heuristic baselines in most NRS literature lag behind the SOTA. Accordingly, baseline heuristics in this experiment are selected from the SOTA heuristic literature and appropriately configured to enable a more informative comparison.

\subsubsection{Does Deep Learning Truly Help in NRSs?}

The experimental results suggest that DL does contribute to good NRS performance. Under the conventional evaluation pipeline, most NRSs explicitly designed for generalization outperform their handcrafted heuristic counterparts within the same category and framework on uniformly distributed instances on different scales. Notably, across both pipelines, appending methods like ReLD, ICAM, and L2R, which incorporate distance information as an auxiliary bias, can achieve stronger generalization performance than the nearest neighbor heuristic. This outperformance suggests that DL models can extract useful implicit knowledge complementary to explicit distance information, demonstrating generalization potential.

On the other hand, the current generalization capability of NRSs remains limited. Under the proposed pipeline, all NRSs demonstrate lower-than-expected effectiveness and, at times, worse reliability. In several cases, their effectiveness even falls below that of simple construction-based heuristics, suggesting that rules learned solely from instances with the uniform node distribution fail to transfer robustly across diverse distributions.

Taken together, these results support a cautiously optimistic conclusion. DL can indeed capture implicit knowledge and yield measurable gains in solving various routing problems. The less favorable performance under the proposed pipeline likely stems from overfitting due to a narrow training distribution, rather than a fundamental limitation of NRSs. In addition, the algorithmic frameworks of current NRSs, especially improvement-based ones, are often simpler than those of advanced heuristics, which may also constrain performance. Therefore, NRSs retain clear research value and foreseeable potential for further performance gain in in-problem generalization. Within the same heuristic framework, DL holds the promise of discovering rules that outperform or complement handcrafted designs, thereby further improving overall performance. Moreover, since few-shot adaptation can rapidly align learned implicit knowledge with a target distribution, NRSs offer a viable pathway to practical deployment by combining general-purpose knowledge with distribution-specific patterns to serve a wide range of applications.

\section{Challenges, Frontier Strategies, and \\ Future Directions}
\label{Section 7: Challenges}

The increasing demand for NRSs to perform effectively and reliably in real-world settings has brought significant challenges in generalization. This section elaborates on these challenges, discusses the strategies explored in the recent literature, and outlines potential future research directions.

\subsection{In-problem Generalization}
\label{Section 7.1: In-problem}

In-problem generalization refers to maintaining stable performance across instances from different data distributions of a single problem, including but not limited to variations in scale, node distribution, and constraint tightness. It is influenced by both how models extract and utilize information and by the distribution shift between the training and test data. Correspondingly, related studies are analyzed from two complementary perspectives: model design and data distribution.

\textbf{Model Design } Strategies for improving in-problem generalization through model design follow two main lines. On the one hand, a few appending methods~\cite{peng2020deep,xin2021multi,drakulic2024bq,luo2023neural,luo2025boosting,fang2024invit,zhou2025l2r} and restricted direct LNS methods~\cite{li2025destroy} allocate more attention layers to dynamically capture relationships between the partial solution and remaining nodes, and among the remaining nodes themselves. On the other hand, some appending methods~\cite{jin2023pointerformer,son2023meta,li2023learning,gao2024towards,wang2024distance,zhou2024instance} incorporate node-wise distances into specific model modules, given VRPs' typical distance-based objectives. Both strategies enable more informed node selection, albeit with different trade-offs. The former incurs higher memory and runtime overhead due to stepwise re-embedding, while the latter relies more on handcrafted designs and is applicable only to distance-based objectives. These limitations highlight the need for future architectures that jointly improve state representation and inference efficiency without restricting to specific objectives.

\textbf{Data Distribution } Existing strategies related to data distribution generally employ two strategies. The first strategy pre-processes test instances to resemble the training distribution. For example, certain two-stage methods~\cite{hou2023generalize,pan2023htsp}, restricted direct LNS methods~\cite{luo2023neural,luo2025boosting,li2021learning,kim2021learning,zong2022rbg,cheng2023select,ye2024glop,zheng2024udc}, and sequential search methods~\cite{fu2021generalize} decompose the original problem into subproblems with scales comparable to those of the training data. Besides, statically~\cite{qiu2022dimes,sun2023difusco,li2024distribution,li2024fast,lischka2024less,fu2021generalize,xing2020solve,min2024unsupervised} or dynamically~\cite{gao2024towards,wang2024distance,fang2024invit,drakulic2024bq} reducing the search space~\cite{sun2021using} has been proven effective across categories. Coordinate normalization can further align the test distribution with the training distributions~\cite{ye2024glop,zheng2024udc,fu2021generalize,fang2024invit,zhou2025l2r,chen2025improving}. However, they often prioritize locally optimal partial solutions, which can trap the complete solution in local optima and degrade final performance. Incorporating global information during search may help mitigate this issue. The second strategy diversifies the training data, \eg, by incorporating instances with different scales or node distributions~\cite{zhou2024instance,wang2024distance,zhou2023towards,khalil2017learning,jiang2022learning,huang2025rethinking,wang2024asp}. Given that NRSs are typically trained on narrow distributions and may overfit, as our experiments suggest, enriching the training data with diverse distributions may already yield substantial gains. Nevertheless, identifying representative distributions and designing effective training strategies remain challenging. Neither the factors shaping data distributions nor the mechanisms through which distributions affect model performance are yet systematically understood. A thorough investigation of these issues is therefore needed in future work.

\subsection{Cross-problem Generalization}

Most NRSs train a specialized model for each problem. This ``one problem, one model'' paradigm is inefficient because it ignores the structural similarities across VRPs. Therefore, a promising research direction is to develop a general-purpose solver that can handle multiple VRPs without costly problem-specific engineering or retraining from scratch.
Existing strategies primarily adapt established DL techniques to the AR appending methods. For example, a recent trend for developing general VRP solvers is multi-task learning~\cite{liu2024multi,zhou2024mvmoe,li2024cada,berto2024routefinder,liu2025mixed,goh2025shield,zheng2025mtlkd}, where a single model is trained and tested on a set of VRP variants with combinations of predefined attributes. However, the observed generalization is at best limited to variants with novel attribute combinations and, at worst, amounts to in-domain performance. This reliance on a predefined attribute set fundamentally restricts generalization, since the attributes in real-world problems cannot be fully anticipated or enumerated in advance. Other strategies employ the model with a shared backbone and problem-specific adapters for different VRPs~\cite{lin2024cross,drakulic2024goal}. While this design enables flexible fine-tuning, it introduces unavoidable limitations. Particularly, related constraint handling remains inherently problem-specific and tied to the adapter design, preventing zero-shot application to unseen problems. Furthermore, this design incurs additional computational overhead during fine-tuning.

To outline potential pathways toward zero-shot generalization for unseen VRPs, two promising future research directions are highlighted: input representations and constraint handling. (1) For \textbf{input representations}, existing methods largely rely on fixed-length attribute or problem vectors, which inherently limit the range of solvable problems. Accordingly, moving beyond attribute-predefined designs toward more general input representations is thus a critical step for broader applicability. (2) For \textbf{constraint handling}, step-wise masking in AR single-stage methods can enforce hard constraints but introduce manual intervention and is inapplicable to certain problems~\cite{bi2024learning} (\eg, TSPTW). To address both issues, a promising direction is to develop intervention-free constraint-handling mechanisms, especially for cross-problem settings. 

\section{Conclusions}
\label{Section 8: Conclusions}

This survey systematically reviews \underline{n}eural \underline{r}outing \underline{s}olvers (NRSs) from the perspective of heuristics. They are identified as heuristic algorithms in which DL-learned rules replace handcrafted ones. A hierarchical taxonomy is introduced based on how solutions are constructed or improved. This perspective enables consistent analysis of connections and developmental trends among NRSs, and naturally links their designs to established heuristic principles within corresponding categories.

Besides, a generalization-focused evaluation pipeline is proposed to address limitations of the conventional one, and representative NRSs are benchmarked under both pipelines. Results under the new pipeline show that NRSs trained on a narrow range of instance distributions can be outperformed by simple construction-based heuristics such as nearest neighbor and random insertion, indicating that the conventional pipeline can lead to overly optimistic conclusions. These findings motivate further discussion of the new pipeline's advantages, principles for method selection, and the role of DL in NRSs despite current performance gaps. Finally, two central challenges in the field, \ie, in-problem and cross-problem generalization, are analyzed. Related prevailing strategies are summarized, and several directions for future work are outlined.

\section*{Acknowledgments}

The authors would like to thank K. Li, Z. Xie, F. Luo, and J. Hou for valuable support in the experimental studies.

\bibliography{ref}
\bibliographystyle{IEEEtran}

% Appendix
\clearpage
\newpage
\appendices

This is the supplementary material for ``Survey on Neural Routing Solvers''. 

\section{Experiment Details}
\label{Appendix: Exp_Details}

\subsection{Adopted Resources}

The sources and possible licenses of the adopted methods and benchmark instances are summarized in Table~\ref{tab:source_method} and~\ref{tab:source_benchmarks}. All of them are open-sourced and available for academic use.

\begin{table*}[htbp]
  \centering
  \caption{Sources of the Adopted Methods}
  \tabcolsep=0.53cm
  % \vspace{-0.2cm}
    \begin{tabular}{lll}
    \toprule[0.5mm]
    \textbf{Method}  & \multicolumn{1}{l}{\textbf{Link}} & \multicolumn{1}{l}{\textbf{License}} \\
    \midrule
    LKH-3~\cite{helsgaun2017extension}    & \url{http://webhotel4.ruc.dk/~keld/research/LKH-3/} & Available for academic research use \\
    HGS~\cite{vidal2022hybrid}      & \url{https://github.com/vidalt/HGS-CVRP} & MIT License \\
    AILS-II~\cite{maximo2024ails}    &    \url{https://github.com/INFORMSJoC/2023.0106}   & MIT License \\
    \midrule
    BQ~\cite{drakulic2024bq}       & \url{https://github.com/naver/bq-nco} & CC BY-NC-SA 4.0 license \\
    LEHD~\cite{luo2023neural}    &   \url{https://github.com/CIAM-Group/NCO_code/tree/main/single_objective/LEHD}  &  MIT License \\
    SIL~\cite{luo2025boosting}    &    \url{https://github.com/CIAM-Group/SIL}   & MIT License \\
    ICAM~\cite{zhou2024instance}     & \url{https://github.com/CIAM-Group/ICAM} & MIT License \\
    ELG~\cite{gao2024towards}      &  \url{https://github.com/gaocrr/ELG} & MIT License \\
    INViT~\cite{fang2024invit}    &  \url{https://github.com/Kasumigaoka-Utaha/INViT}  & MIT License  \\
    L2R~\cite{zhou2025l2r}      &   \url{https://github.com/CIAM-Group/L2R}    & MIT License \\
    DGL~\cite{xiao2025dgl}      &   \url{https://github.com/wuyuesong/DGL}    & Available for academic research use \\
    ReLD~\cite{huang2025rethinking}     &   \url{https://github.com/ziweileonhuang/reld-nco}  & MIT License \\
    L2C-Insert~\cite{luo2025learning}    &    \url{https://github.com/CIAM-Group/L2C_Insert}   & MIT License \\
    H-TSP~\cite{pan2023htsp}    &    \url{https://github.com/Learning4Optimization-HUST/H-TSP}   & MIT License \\
    DACT~\cite{ma2021learning}     &    \url{https://github.com/yining043/VRP-DACT}   & MIT License \\
    NeuOpt~\cite{ma2024learning}    &   \url{https://github.com/yining043/NeuOpt}    & MIT License \\
    GenSCO~\cite{li2025generation}    &  \url{https://github.com/Thinklab-SJTU/GenSCO}    & Available for academic research use \\
    DRHG~\cite{li2025destroy}     &   \url{https://github.com/CIAM-Group/DRHG}   & Available for academic research use \\
    FastT2T~\cite{li2024fast}    &  \url{https://github.com/Thinklab-SJTU/Fast-T2T}   & MIT license \\
    GFACS~\cite{kim2024ant}    &    \url{https://github.com/ai4co/gfacs}   & MIT license \\
    \bottomrule[0.5mm]
    \end{tabular}%
  \vspace{-0.3cm}
  \label{tab:source_method}%
\end{table*}%

\begin{table*}[htbp]
  \centering
  \caption{Sources of the Adopted Benchmarks}
  \tabcolsep=0.36cm
  % \vspace{-0.2cm}
    \begin{tabular}{p{7em}p{22em}lll}
    \toprule[0.5mm]
    \multicolumn{1}{l}{\textbf{Benchmark}} & \multicolumn{1}{l}{\textbf{Instance}} & \multicolumn{1}{l}{\textbf{BKS}} \\
    \midrule
    TSPLIB &   \url{http://comopt.ifi.uni-heidelberg.de/software/TSPLIB95/}    & \url{http://comopt.ifi.uni-heidelberg.de/software/TSPLIB95/STSP.html} \\
    \midrule
    National &     \url{https://www.math.uwaterloo.ca/tsp/world/countries.html}  & \url{https://www.math.uwaterloo.ca/tsp/world/summary.html} \\
    \midrule
    VLSI  &    \url{https://www.math.uwaterloo.ca/tsp/vlsi/index.html}   & \url{https://www.math.uwaterloo.ca/tsp/vlsi/summary.html} \\
    \midrule
    8th DIMACS & \multirow{2}{*}{\makecell[{{p{22em}}}]{\url{http://dimacs.rutgers.edu/archive/Challenges/TSP/download.html}}} & \url{http://webhotel4.ruc.dk/~keld/research/LKH/DIMACS_results.html} \\
    Implementation Challenge &      & \url{http://dimacs.rutgers.edu/archive/Challenges/TSP/opts.html} \\
    \midrule
    CVRPLIB &    \url{https://galgos.inf.puc-rio.br/cvrplib/index.php/en/instances}   & Provided in the corresponding .vrp files of the instances. \\
    \bottomrule[0.5mm]
    \end{tabular}%
    % \vspace{-0.3cm}
  \label{tab:source_benchmarks}%
\end{table*}%

\subsection{Detailed Results of the Proposed Pipeline} The detailed results of the proposed pipeline, including NRSs with all available models, are presented in Table~\ref{tab:exp_proposed_TSP_details} and Table~\ref{tab:exp_proposed_CVRP_details}, respectively.

\begin{table*}[t]
  \centering
  \caption{Detailed Experimental Results of the Proposed Evaluation Pipeline for TSP}
  % \vspace{-0.2cm}
  \tabcolsep=0.185cm
  \begin{threeparttable}
    \begin{tabular}{l|ccc|ccc|ccc|cc}
    \toprule[0.5mm]
    \multicolumn{1}{c|}{\multirow{2}[2]{*}{Method}} & \multicolumn{3}{c|}{(0,1K)} & \multicolumn{3}{c|}{[1K, 10K)} & \multicolumn{3}{c|}{[10K, 100K]} & \multicolumn{2}{c}{Total} \\
\cmidrule{2-12}          & Gap   & Time  & \multicolumn{1}{l|}{Solved} & Gap   & Time  & \multicolumn{1}{l|}{Solved} & Gap   & Time  & \multicolumn{1}{l|}{Solved} & Gap   & Solved \\
    \midrule
    Nearest Neighbor & 25.29\% & 0.01s & 69/69 & 26.66\% & 0.29s & 109/109 & 25.01\% & 22.60s & 50/50 & 25.88\% & 228/228 \\
    Random Insertion & 10.60\% & 0.00s & 69/69 & 15.32\% & 0.05s & 109/109 & 16.37\% & 8.93s & 50/50 & 14.12\% & 228/228 \\
    LKH-3\phantom{$^{\downarrow}$} t=n/3, runs=1 & 0.00\% & 7.88s & 69/69 & 0.01\% & 631.34s & 109/109 & 0.08\% & 14600.50s & 50/50 & 0.03\% & 228/228 \\
    LKH-3$^\downarrow$ t=n/3, runs=1 & 0.00\% & 9.25s & 69/69 & 0.01\% & 600.36s & 109/109 & 0.05\% & 10800.24s & 50/50 & 0.02\% & 228/228 \\
    \midrule
    BQ    & 5.00\% & 2.51s & 68/69 & 19.03\% & 22.74s & 92/109 & 52.00\% & 187.81s & 4/50  & 14.02\% & 164/228 \\
    LEHD$^*$ greedy & 4.85\% & 1.01s & 69/69 & 20.13\% & 68.27s & 106/109 & 49.35\% & 1386.01s & 11/50 & 16.19\% & 186/228 \\
    SIL$^*$ greedy (1K) & 6.07\% & 1.72s & 69/69 & 10.29\% & 17.31s & 106/109 & 20.71\% & 434.58s & 47/50 & 11.18\% & 222/228 \\
    SIL$^*$ greedy (5K) & 9.68\% & 1.71s & 69/69 & 11.06\% & 17.83s & 109/109 & 18.75\% & 422.26s & 50/50 & 12.33\% & 228/228 \\
    SIL$^*$ greedy (10K) & 6.11\% & 1.72s & 69/69 & 9.99\% & 17.87s & 109/109 & 20.85\% & 414.78s & 50/50 & 11.20\% & 228/228 \\
    SIL$^*$ greedy (50K) & 9.86\% & 1.71s & 69/69 & 9.40\% & 17.87s & 109/109 & 12.48\% & 427.04s & 50/50 & 10.21\% & 228/228 \\
    SIL$^*$ greedy (100K) & 8.64\% & 1.69s & 69/69 & 9.83\% & 17.69s & 109/109 & 11.11\% & 430.73s & 50/50 & 9.75\% & 228/228 \\
    ICAM  & 6.53\% & 0.25s & 69/69 & 16.62\% & 21.57s & 109/109 & 21.34\% & 1050.33s & 19/50 & 13.54\% & 197/228 \\
    ELG   & 6.05\% & 0.63s & 69/69 & 18.14\% & 88.12s & 108/109 & 21.65\% & 940.35s & 6/50  & 13.70\% & 183/228 \\
    INViT-3V   & 7.93\% & 2.77s & 69/69 & 12.08\% & 49.03s & 109/109 & 11.52\% & 1079.50s & 42/50 & 10.67\% & 220/228 \\
    L2R   & 5.89\% & 1.60s & 69/69 & 9.22\% & 15.55s & 109/109 & 8.52\% & 153.11s & 50/50 & 8.06\% & 228/228 \\
    DGL   & 6.53\% & 1.17s & 69/69 & 11.32\% & 11.67s & 109/109 & 11.14\% & 58.62s & 25/50 & 9.67\% & 203/228 \\
    \midrule
    L2C-Insert$^*$ greedy & 4.39\% & 1.51s & 69/69 & 18.12\% & 15.34s & 109/109 & 30.94\% & 145.17s & 50/50 & 16.77\% & 228/228 \\
    \midrule
    H-TSP (1K) & 6.16\% & 0.67s & 36/69 & 12.26\% & 3.24s & 100/109 & 12.74\% & 21.93s & 40/50 & 11.12\% & 176/228 \\
    H-TSP (2K) & 6.26\% & 0.64s & 36/69 & 11.83\% & 3.23s & 100/109 & 12.35\% & 21.42s & 40/50 & 10.81\% & 176/228 \\
    H-TSP (5K) & 6.16\% & 0.61s & 36/69 & 11.62\% & 3.15s & 100/109 & 12.29\% & 21.44s & 40/50 & 10.65\% & 176/228 \\
    H-TSP (10K) & 6.10\% & 0.61s & 36/69 & 11.75\% & 3.13s & 100/109 & 12.29\% & 21.37s & 40/50 & 10.72\% & 176/228 \\
    \midrule
    DACT T=1K (20) & 24.54\% & 39.47s & 69/69 & 26.85\% & 260.55s & 83/109 & /     & /     & 0/50  & 25.80\% & 152/228 \\
    DACT T=1K (50) & 17.49\% & 39.68s & 69/69 & 26.69\% & 259.77s & 83/109 & /     & /     & 0/50  & 22.52\% & 152/228 \\
    DACT T=1K (100) & 16.37\% & 39.84s & 69/69 & 26.58\% & 261.73s & 83/109 & /     & /     & 0/50  & 21.94\% & 152/228 \\
    \midrule
    NeuOpt T=1K (20) & 35.76\% & 533.78s & 7/69  & /     & /     & 0/109 & /     & /     & 0/50  & 35.76\% & 7/228 \\
    NeuOpt T=1K (50) & 20.70\% & 138.26s & 27/69 & /     & /     & 0/109 & /     & /     & 0/50  & 20.70\% & 27/228 \\
    NeuOpt T=1K (100) & 12.44\% & 104.05s & 36/69 & /     & /     & 0/109 & /     & /     & 0/50  & 12.44\% & 36/228 \\
    NeuOpt T=1K (200) & 19.90\% & 81.22s & 46/69 & /     & /     & 0/109 & /     & /     & 0/50  & 19.90\% & 46/228 \\
    \midrule
    L2C-Insert$^*$ T=1K & 1.08\% & 381.55s & 69/69 & 9.80\% & 479.98s & 109/109 & 29.15\% & 615.19s & 50/50 & 11.41\% & 228/228 \\
    GenSCO (100) & 14.56\% & 23.19s & 68/69 & 35.46\% & 677.31s & 104/109 & 35.17\% & 14304.70s & 25/50 & 28.21\% & 197/228 \\
    GenSCO (500) & 19.64\% & 25.57s & 42/69 & 48.66\% & 563.68s & 68/109 & 52.89\% & 17964.48s & 14/50 & 39.31\% & 124/228 \\
    GenSCO (1K) & 24.11\% & 27.76s & 29/69 & 19.27\% & 268.60s & 53/109 & 90.42\% & 24542.75s & 2/50  & 22.64\% & 84/228 \\
    \midrule
    LEHD$^*$ RRC1K & 1.73\% & 498.40s & 69/69 & 10.87\% & 1634.51s & 109/109 & 24.02\% & 2769.86s & 9/50  & 8.13\% & 187/228 \\
    SIL$^*$ PRC1K (1K) & 0.63\% & 903.72s & 69/69 & 2.89\% & 3250.51s & 109/109 & 7.29\% & 5366.67s & 50/50 & 3.17\% & 228/228 \\
    SIL$^*$ PRC1K (5K) & 0.86\% & 888.48s & 69/69 & 2.86\% & 2889.38s & 109/109 & 6.85\% & 5046.36s & 50/50 & 3.13\% & 228/228 \\
    SIL$^*$ PRC1K (10K) & 0.55\% & 891.35s & 69/69 & 2.63\% & 2898.09s & 109/109 & 6.36\% & 4962.70s & 50/50 & 2.82\% & 228/228 \\
    SIL$^*$ PRC1K (50K) & 0.80\% & 883.45s & 69/69 & 2.56\% & 2879.65s & 109/109 & 5.11\% & 5009.26s & 50/50 & 2.59\% & 228/228 \\
    SIL$^*$ PRC1K (100K) & 0.80\% & 883.87s & 69/69 & 2.58\% & 2880.46s & 109/109 & 4.55\% & 4933.45s & 50/50 & 2.47\% & 228/228 \\
    DRHG T=1K & 0.10\% & 769.53s & 69/69 & 1.46\% & 2857.55s & 109/109 & 4.46\% & 3004.28s & 50/50 & 1.71\% & 228/228 \\
\midrule    Fast T2T T\textsubscript{s}=10, T\textsubscript{g}=10 (50) & 15.72\% & 1.42s & 45/69 & /     & /     & 0/109 & /     & /     & 0/50  & 15.72\% & 45/228 \\
    Fast T2T T\textsubscript{s}=10, T\textsubscript{g}=10 (100) & 10.46\% & 1.41s & 45/69 & /     & /     & 0/109 & /     & /     & 0/50  & 10.46\% & 45/228 \\
    Fast T2T T\textsubscript{s}=10, T\textsubscript{g}=10 (500) & 20.68\% & 1.39s & 45/69 & /     & /     & 0/109 & /     & /     & 0/50  & 20.68\% & 45/228 \\
    Fast T2T T\textsubscript{s}=10, T\textsubscript{g}=10 (1K) & 18.55\% & 1.40s & 45/69 & /     & /     & 0/109 & /     & /     & 0/50  & 18.55\% & 45/228 \\
    Fast T2T T\textsubscript{s}=10, T\textsubscript{g}=10 (10K) & 43.74\% & 1.41s & 45/69 & /     & /     & 0/109 & /     & /     & 0/50  & 43.74\% & 45/228 \\
    \midrule
    GFACS$^{\dagger}$ T=100, K=100 (200) & 31.64\% & 166.94s & 66/69 & 86.77\% & 2601.13s & 22/109 & /     & /     & 0/50  & 45.42\% & 88/228 \\
    GFACS$^{\dagger}$ T=100, K=100 (500) & 35.06\% & 175.46s & 64/69 & 86.43\% & 2590.34s & 16/109 & /     & /     & 0/50  & 45.33\% & 80/228 \\
    GFACS$^{\dagger}$ T=100, K=100 (1K) & 42.05\% & 175.18s & 64/69 & 84.88\% & 2766.93s & 10/109 & /     & /     & 0/50  & 47.84\% & 74/228 \\
    GFACS$^{\ddagger}$ T=100, K=100 (200) & 0.72\% & 174.21s & 69/69 & 3.76\% & 9142.93s & 83/109 & /     & /     & 0/50  & 2.38\% & 152/228 \\
    GFACS$^{\ddagger}$ T=100, K=100 (500) & 0.80\% & 209.26s & 69/69 & 3.31\% & 6718.88s & 65/109 & /     & /     & 0/50  & 2.02\% & 134/228 \\
    GFACS$^{\ddagger}$ T=100, K=100 (1K) & 0.93\% & 225.39s & 69/69 & 3.63\% & 7545.88s & 66/109 & /     & /     & 0/50  & 2.25\% & 135/228 \\
    \bottomrule[0.5mm]
    \end{tabular}%
    \begin{tablenotes}
      \footnotesize
      \item[]%
        \hspace*{-\dimexpr\labelwidth+\labelsep\relax}%
        \makebox[\labelwidth][l]{()}%
        \hspace{\labelsep}%
        \hspace*{0.2em}The content in parentheses indicates the node scale of instances with uniformly distributed nodes in the model's training set.
      \item[$\downarrow$] \hspace*{0.2em}$Initial\_Period$ is set as 1K, rather than the original value DIMENSION/2.
      \item[$*$]  \hspace*{0.2em}The NRS supports more than one inference strategy.
      \item[$\dagger$] \hspace*{0.2em}GFACS without local search at the last generation. The output solution is the best of the population at the last generation.
      \item[$\ddagger$] \hspace*{0.2em}The original version of GFACS. The output solution is the best in history (all after local search).
    \end{tablenotes}
    \end{threeparttable}
    % \vspace{-0.3cm}
  \label{tab:exp_proposed_TSP_details}%
\end{table*}%

\begin{table*}[t]
  \centering
  \caption{Detailed Experimental Results of the Proposed Evaluation Pipeline for CVRP}
  % \vspace{-0.2cm}
  \tabcolsep=0.18cm
  \begin{threeparttable}
    \begin{tabular}{l|ccc|ccc|ccc|cc}
    \toprule[0.5mm]
    \multicolumn{1}{c|}{\multirow{2}[2]{*}{Method}} & \multicolumn{3}{c|}{(0,1K)} & \multicolumn{3}{c|}{[1K, 10K)} & \multicolumn{3}{c|}{[10K, 100K]} & \multicolumn{2}{c}{Total} \\
\cmidrule{2-12}          & Gap   & Time  & Solved & Gap   & Time  & Solved & Gap   & Time  & \multicolumn{1}{l|}{Solved} & Gap   & Solved \\
    \midrule
    Nearest Neighbor & 21.17\% & 0.03s & 99/99 & 15.18\% & 1.08s & 5/5   & 11.80\% & 14.63s & 6/6   & 20.39\% & 110/110 \\
    Random Insertion & 75.00\% & 0.00s & 36/99 & /     & /     & 0/5   & /     & /     & 0/6   & 75.00\% & 36/110 \\    HGS t=n/3 & 0.29\% & 111.24s & 99/99 & 3.59\% & 1428.41s & 5/5   & 7.86\% & 6926.41s & 6/6   & 0.85\% & 110/110 \\
    AILS-II t=n/3 & 0.57\% & 133.95s & 99/99 & 1.58\% & 1388.42s & 5/5   & 1.58\% & 5646.48s & 6/6   & 0.68\% & 110/110 \\
    \midrule
    BQ    & 8.87\% & 3.63s & 99/99 & 20.28\% & 39.97s & 5/5   & 41.52\% & 202.69s & 5/6   & 10.89\% & 109/110 \\
    {LEHD$^*$ greedy} & 11.25\% & 1.53s & 98/99 & 19.22\% & 99.43s & 5/5   & 32.80\% & 852.02s & 2/6   & 12.04\% & 105/110 \\
    {SIL$^*$ greedy (1K)} & 43.24\% & 2.62s & 58/99 & 21.73\% & 27.17s & 5/5   & 15.39\% & 146.85s & 6/6   & 39.26\% & 69/110 \\
    {SIL$^*$ greedy (5K)} & 42.70\% & 2.49s & 62/99 & 18.34\% & 26.60s & 5/5   & 12.38\% & 147.91s & 6/6   & 38.54\% & 73/110 \\
    {SIL$^*$ greedy (10K)} & 46.94\% & 2.54s & 62/99 & 20.78\% & 27.08s & 5/5   & 12.62\% & 148.55s & 6/6   & 42.32\% & 73/110 \\
    {SIL$^*$ greedy (50K)} & 40.04\% & 2.48s & 65/99 & 16.09\% & 26.86s & 5/5   & 10.81\% & 146.83s & 6/6   & 36.16\% & 76/110 \\
    {SIL$^*$ greedy (100K)} & 40.44\% & 2.54s & 60/99 & 18.32\% & 26.86s & 5/5   & 10.95\% & 148.69s & 6/6   & 36.39\% & 71/110 \\
    ICAM  & 5.00\% & 0.42s & 99/99 & 11.69\% & 32.32s & 5/5   & /     & /     & 0/6   & 5.32\% & 104/110 \\
    ELG   & 8.03\% & 1.29s & 99/99 & 18.51\% & 30.21s & 5/5   & 29.38\% & 133.08s & 2/6   & 8.93\% & 106/110 \\
    INViT-3V   & 13.15\% & 4.72s & 99/99 & 19.03\% & 77.33s & 5/5   & 23.91\% & 496.25s & 5/6   & 13.91\% & 109/110 \\
    L2R   & 8.16\% & 2.49s & 99/99 & 11.62\% & 23.99s & 5/5   & 11.08\% & 97.12s & 6/6   & 8.48\% & 110/110 \\
    DGL   & 15.27\% & 2.22s & 99/99 & 17.96\% & 22.60s & 5/5   & 18.69\% & 78.64s & 5/6   & 15.55\% & 109/110 \\
    ReLD  & 4.10\% & 0.41s & 99/99 & 10.22\% & 5.29s & 5/5   & 11.27\% & 28.87s & 3/6   & 4.58\% & 107/110 \\
    \midrule
    {L2C-Insert$^*$ greedy} & 6.87\% & 2.73s & 99/99 & 22.37\% & 616.75s & 5/5   & 49.41\% & 5525.71s & 2/6   & 8.40\% & 106/110 \\
    \midrule
    DACT T=1K (20) & 18.90\% & 284.37s & 64/99 & /     & /     & 0/5   & /     & /     & 0/6   & 18.90\% & 64/110 \\
    DACT T=1K (50) & 16.42\% & 338.69s & 54/99 & /     & /     & 0/5   & /     & /     & 0/6   & 16.42\% & 54/110 \\
    DACT T=1K (100) & 16.42\% & 246.51s & 74/99 & 17.70\% & 479.82s & 1/5   & /     & /     & 0/6   & 16.44\% & 75/110 \\
    \midrule
    NeuOpt T=1K (20) & 79.48\% & 3440.62s & 6/99  & /     & /     & 0/5   & /     & /     & 0/6   & 79.48\% & 6/110 \\
    NeuOpt T=1K (50) & 74.85\% & 2953.52s & 7/99  & /     & /     & 0/5   & /     & /     & 0/6   & 74.85\% & 7/110 \\
    NeuOpt T=1K (100) & 46.93\% & 868.55s & 24/99 & /     & /     & 0/5   & /     & /     & 0/6   & 46.93\% & 24/110 \\
    NeuOpt T=1K (200) & 26.93\% & 571.14s & 36/99 & /     & /     & 0/5   & /     & /     & 0/6   & 26.93\% & 36/110 \\
    \midrule
    L2C-Insert$^*$ T=1K & 3.21\% & 344.05s & 99/99 & 18.87\% & 6166.21s & 5/5   & 44.29\% & 32754.86s & 2/6   & 4.72\% & 106/110 \\
    \midrule
    {LEHD$^*$ RRC1K} & 3.58\% & 796.15s & 99/99 & 11.73\% & 2043.74s & 5/5   & 21.98\% & 2820.28s & 2/6   & 4.32\% & 106/110 \\
    SIL$^*$ PRC1K (1K) & 23.67\% & 1291.37s & 99/99 & 12.59\% & 3405.06s & 5/5   & 11.46\% & 4254.08s & 6/6   & 22.50\% & 110/110 \\
    {SIL$^*$ PRC1K (5K)} & 21.58\% & 1314.59s & 99/99 & 10.12\% & 3487.80s & 5/5   & 8.60\% & 4324.59s & 6/6   & 20.35\% & 110/110 \\
    {SIL$^*$ PRC1K (10K)} & 22.16\% & 1297.48s & 99/99 & 9.40\% & 3442.22s & 5/5   & 8.54\% & 4223.01s & 6/6   & 20.84\% & 110/110 \\
    SIL$^*$ PRC1K (50K) & 21.38\% & 1307.97s & 99/99 & 8.28\% & 3471.69s & 5/5   & 7.40\% & 4251.88s & 6/6   & 20.02\% & 110/110 \\
    {SIL$^*$ PRC1K (100K)} & 22.24\% & 1315.62s & 99/99 & 9.13\% & 3471.46s & 5/5   & 8.02\% & 4267.11s & 6/6   & 20.87\% & 110/110 \\
    DRHG T=1K & 11.11\% & 1114.60s & 99/99 & 17.95\% & 2529.12s & 5/5   & 16.95\% & 5376.37s & 6/6   & 11.74\% & 110/110 \\
    \midrule
    GFACS$^{\dagger}$ T=100, K=100 (200) & 36.83\% & 437.38s & 99/99 & 34.09\% & 9654.33s & 3/5   & /     & /     & 0/6   & 36.75\% & 102/110 \\
    GFACS$^{\dagger}$ T=100, K=100 (500) & 61.11\% & 606.05s & 78/99 & 86.40\% & 1701.41s & 1/5   & /     & /     & 0/6   & 61.43\% & 79/110 \\
    GFACS$^{\dagger}$ T=100, K=100 (1K) & 67.68\% & 476.95s & 64/99 & /     & /     & 0/5   & /     & /     & 0/6   & 67.68\% & 64/110 \\
    GFACS$^{\ddagger}$ T=100, K=100 (200) & 2.60\% & 405.81s & 99/99 & 7.65\% & 14884.94s & 4/5   & /     & /     & 0/6   & 2.80\% & 103/110 \\
    GFACS$^{\ddagger}$ T=100, K=100 (500) & 2.48\% & 588.86s & 99/99 & 5.28\% & 10935.70s & 3/5   & /     & /     & 0/6   & 2.56\% & 102/110 \\
    GFACS$^{\ddagger}$ T=100, K=100 (1K) & 2.55\% & 511.39s & 99/99 & 5.50\% & 11283.58s & 3/5   & /     & /     & 0/6   & 2.63\% & 102/110 \\
    \bottomrule[0.5mm]
    \end{tabular}%
    \begin{tablenotes}
      \footnotesize
      \item[]%
        \hspace*{-\dimexpr\labelwidth+\labelsep\relax}%
        \makebox[\labelwidth][l]{()}%
        \hspace{\labelsep}%
        \hspace*{0.2em}The content in parentheses indicates the node scale of instances with uniformly distributed nodes in the model's training set.
      \item[$*$]  \hspace*{0.2em}The NRS supports more than one inference strategy.
      \item[$\dagger$] \hspace*{0.2em}GFACS without local search at the last generation. The output solution is the best of the population at the last generation.
      \item[$\ddagger$] \hspace*{0.2em}The original version of GFACS. The output solution is the best in history (all after local search).
    \end{tablenotes}
    \end{threeparttable}
    % \vspace{-0.3cm}
  \label{tab:exp_proposed_CVRP_details}%
\end{table*}%

\vfill
\end{CJK*}
\end{document}